\documentclass[12pt]{article}
\usepackage{amsmath}
\usepackage{latexsym}
\usepackage{amssymb}
\usepackage{a4wide}
\newtheorem{thm}{Theorem}[section]
\newtheorem{la}[thm]{Lemma}
\newtheorem{Defn}[thm]{Definition}
\newtheorem{Remark}[thm]{Remark}
\newtheorem{Note}[thm]{Note}
\newtheorem{prop}[thm]{Proposition}

\newtheorem{Example}[thm]{Example}
\newtheorem{Examples}[thm]{Examples}
\newtheorem{Problems}[thm]{Problems}

\newtheorem{Problem}[thm]{Problem}
\newtheorem{Number}[thm]{\!\!}
\newenvironment{defn}{\begin{Defn}\rm}{\end{Defn}}

\newenvironment{rem}{\begin{Remark}\rm}{\end{Remark}}
\newenvironment{numba}{\begin{Number}\rm}{\end{Number}}
\newenvironment{proof}{{\noindent\bf Proof.}}%
                  {\nopagebreak\hspace*{\fill}$\Box$\medskip\medskip\par}   
\newcommand{\Punkt}{\nopagebreak\hspace*{\fill}$\Box$}

\newcommand{\wt}{\widetilde}

\newcommand{\n}{\rm}

\newcommand{\mto}{\mapsto}

\newcommand{\ve}{\varepsilon}

\newcommand{\N}{{\mathbb N}}

\newcommand{\R}{{\mathbb R}}

\newcommand{\Q}{{\mathbb Q}}

\newcommand{\Z}{{\mathbb Z}}
\newcommand{\C}{{\mathbb C}}

\newcommand{\K}{{\mathbb K}}

\newcommand{\cO}{{\cal O}}

\newcommand{\sub}{\subseteq}

\newcommand{\cL}{{\cal L}}

\newcommand{\smin}{\mbox{\n\footnotesize min}}

\newcommand{\sbull}{{\scriptscriptstyle \bullet}}

\newcommand{\half}{{\textstyle \frac{1}{2}}}
\newcommand{\quart}{{\textstyle \frac{1}{4}}}
\begin{document}
\begin{center}
{\Large \bf H\"{o}lder continuous homomorphisms between\vspace{2.2mm}
infinite-dimensional Lie groups are smooth}\vspace{4.5 mm}\\
{\bf Helge Gl\"{o}ckner}\vspace{.9mm}
\end{center}
\noindent{\bf Abstract.\/}
Let
$f \!: G\to H$ be a homomorphism
between smooth Lie groups modelled on
Mackey complete, locally convex real topological
vector spaces.
We show that if $f$
is H\"{o}lder continuous at~$1$,
then $f$ is smooth.\\[3mm]
{\footnotesize
}
\section*{Introduction}
\noindent
While specific examples of infinite-dimensional
Lie groups have been studied extensively and
are well understood,
in the general theory of infinite-dimensional Lie groups
even very fundamental questions are still open.
Various important unsolved problems
were recorded in the preprint~\cite{MiP} by John Milnor in 1982;
most of them
have resisted all attempts at a solution so far.\footnote{Milnor's second
problem
(does every closed subalgebra
correspond to an immersed Lie subgroup\,?)
had in fact already been solved earlier
by Omori~\cite{Omo} (in the negative).
The other three
main problems remain open.
Various smaller problems mentioned
in Milnor's preprint could be settled:
A Lie group whose exponential map is a local diffeomorphism
at $0$ need not be of Campbell-Hausdorff type
\cite[\S3.4.1]{Rob}.
A real or complex analytic Lie group need not be of Campbell-Hausdorff
type \cite[Rem.\,4.7\,(b)]{FUN}.
The complexification of an enlargible real Banach-Lie algebra
need not be enlargible
\cite[Ex.\,VI.4]{GaN}.
A connected Lie group modelled on a
locally convex space is abelian if and only
if its Lie algebra is abelian~\cite[Prop.\,22.15]{INF}.}
In the present article, we give a partial answer
to Milnor's third problem:
``\,Is a continuous homomorphism between Lie groups
necessarily smooth\,?\,''
As our main result, we show
that every {\em H\"{o}lder\/} continuous
homomorphism is smooth.
More precisely:\\[3mm]
{\bf Main Theorem.}
{\em Let $f\!: G\to H$ be a homomorphism
between Lie groups modelled on real locally convex spaces.
If $f$ is H\"{o}lder continuous at~$1$ and
the modelling space of~$H$ is Mackey complete,
then $f$ is smooth.}\\[3mm]
In particular,
the Main Theorem applies to
{\em Lipschitz\/} continuous homomorphisms.
Milnor only considered Lie groups
modelled on complete
locally convex spaces. Mackey completeness
is a very natural and useful
weakened completeness condition \cite{KaM}.\\[3mm]
{\bf A simple special case.}
The basic idea underlying our approach is
most easily explained for
one-parameter groups.
It is helpful to keep this simplest special case
in mind as a guideline also when dealing with the
general case (which is much harder\!).
Thus, consider a continuous homomorphism $\xi\!: \R\to G$
from $\R$ to a Lie group $G$ modelled on a locally convex space~$E$.
Using a chart, we identify an open identity neighbourhood of~$G$
with an\vspace{.4mm}\pagebreak

\noindent
open $0$-neighbourhood in~$E$.
Making use of the first order Taylor expansion
\[
x^2=2x+R(x)
\]
of the squaring map around the identity~$0$,
for small $t\in \R$ we obtain
$\xi(t)=\xi(\half t)^2
=2\xi(\half t)+R(\xi(\half t))$ and thus
$\xi(\half t)=\half \xi(t)-\half R(\xi(\half t))$.
Applying this formula twice yields
\[
\xi(\quart t)\,=\, \half \xi(\half t)-\half R(\xi(\quart t))
\, =\, \quart \xi(t)-\quart R(\xi(\half t))-\half R(\xi(\quart t))\,.
\]
Similarly, $\xi(2^{-n}t)=2^{-n}\xi(t)-\sum_{k=1}^n2^{k-n-1}R(\xi(2^{-k}t))$
for all $n\in \N$, by induction.
After re-parametrizing~$\xi$,
we may assume that $t=1$ can be chosen here. This gives\vspace{-.3mm}
\begin{equation}\label{fm1}
\frac{\xi(2^{-n})}{2^{-n}}\,=\, \xi(1)-\sum_{k=1}^n 2^{k-1}R(\xi(2^{-k}))
\qquad\mbox{for all $n\in \N$.}\vspace{-1.2mm}
\end{equation}
Now assume that $\xi$ is H\"{o}lder
continuous at~$0$, with H\"{o}lder exponent $\alpha\in \;]0,1]$.
Then $\xi(2^{-k})$ is of order $\cO(2^{-k\alpha})$
(as $k\to\infty$). A first order Taylor remainder
being at most quadratic in the order of its argument,
we see that $R(\xi(2^{-k}))$ is of order $\cO(2^{-2k\alpha})$.
Therefore the summands
$2^{k-1}R(\xi(2^{-k}))$ in (\ref{fm1}) are
of order
$\cO(2^{(1-2\alpha)k})$.
If $\alpha\in \;]\half,1]$,
the preceding estimates show
that $n\mto \sum_{k=1}^n2^{k-1}R(\xi(2^{-k}))$\vspace{-.3mm}
is a Mackey-Cauchy sequence in~$E$ and hence convergent
if $E$ is Mackey complete.
Thus $\lim_{n\to \infty}\frac{\xi(2^{-n})}{2^{-n}}$
exists in~$E$, and apparently this limit
gives us a candidate for $\xi'(0)$.
Of course, it remains to show that $\xi'(0)$ really exists
(this is less obvious\,!),
and
that existence of $\xi'(0)$
entails smoothness of $\xi$.
Also, it remains to remove the requirement that $\alpha>\half$
(but all of this can be done).\\[3mm]
{\bf Organization of the paper.}
After a brief description of
the setting of differential calculus
used in the paper,
in Section~\ref{secbase}
we discuss various properties
a mapping between open subsets of locally convex spaces
(or manifolds) can have at a given point:
H\"{o}lder continuity at~$x$,
total differentiability at~$x$,
and feeble differentiability
(an auxiliary notion which we introduce for internal use).
In Section~\ref{sechom}, we show that
$C^1$-homomorphisms between
Lie groups modelled on real locally convex spaces
are smooth (Lemma~\ref{C1smooth}), and we show that a homomorphism is
$C^1$
if it is
totally (or merely feebly)
differentiable at~$1$ (Lemma~\ref{totsmooth}).
Section~\ref{sechold} is devoted to the proof
of the Main Theorem (Theorem~\ref{main}).
In view of the reduction steps already performed,
the crucial point will be to
deduce total differentiability at~$1$
from H\"{o}lder continuity at~$1$.
Although our main result concerns
real Lie groups,
some of our considerations
are not restricted to the real case
and have been formulated
more generally for valued fields
(at little extra cost).
This enables us to show
in Section~\ref{secpad}
that
H\"{o}lder continuous homomorphisms between $p$-adic
Lie groups modelled on Mackey complete
polynormed $\Q_p$-vector spaces are $C^1$ (Theorem~\ref{varmain}).
Proofs for various auxiliary results,
which are best taken on faith
on a first reading, are compiled
in two appendices.\\[3mm]
%
%
{\bf Analogues in convenient differential calculus.}
In the subsequent paper~\cite{CON},
variants of the ideas presented here
are used 
to show that every $\cL ip^0$-homomorphism between
Lie groups in the sense of convenient differential calculus
(as in \cite{KaM}) is smooth in the convenient sense.
More generally, this
conclusion holds for ``conveniently
H\"{o}lder'' homomorphisms~\cite{CON}.
\section{Basic definitions and facts}\label{secbase}
We compile and develop
basic material.
The proofs are recorded in Appendix~\ref{secapp}.
\begin{center}
{\bf Differential calculus in topological vector spaces}
\end{center}
We are working in
the framework of
differential calculus
known as
Keller's $C^\infty_c$-theory~\cite{Kel}
(going back to Michal and Bastiani),
as used in \cite{RES},
\cite{Ham},
\cite{Mic},
\cite{MiP}, \cite{Mil}
and generalized to a differential
calculus over topological fields in~\cite{BGN}.
We recall some of the basic ideas.
\begin{numba}\label{defnMB}
Let $E$ be a real topological vector space,
$F$ be a real locally convex space,
and $U\sub E$ be open.
A map $f\!: U\to F$ is called $C^1$
if it is continuous, the directional derivative
$df(x,y):=\frac{d}{dt}\big|_{t=0}f(x+ty)$
exists for all $x\in U$ and $y\in E$,
and the mapping $df\!: U\times E\to F$ so obtained
is continuous. Inductively, we say that $f$ is $C^{k+1}$
(for $k\geq 1$)
if $f$ is $C^1$ and $df\!: U\times E\to F$ is $C^k$.
The map $f$ is called $C^\infty$ or {\em smooth\/}
if it is $C^k$ for all $k\in \N$.
\end{numba}
\begin{numba}\label{linktoBGN}
If
$f\!: E\supseteq U\to F$ as before
is $C^1$, define
$f^{[1]}\!: U^{[1]}\to F$
on the open set $U^{[1]}:=\{(x,y,t)\in U\times E\times \R\!:
x+ty\in U\}\sub E\times E\times \R$ via
$f^{[1]}(x,y,t)=\frac{1}{t}(f(x+ty)-f(x))$ if $t\not=0$,
$f^{[1]}(x,y,0):=df(x,y)$.
Then $f^{[1]}$ is continuous, because
for small~$t$ we have the integral representation
$f^{[1]}(x,y,t)=\int_0^1df(x+sty,y)\, ds$,
by the Mean Value Theorem.
Furthermore, by definition of $f^{[1]}$,
\begin{equation}\label{deff1}
f^{[1]}(x,y,t)
=\frac{1}{t}(f(x+ty)-f(x))
\qquad\mbox{for all $(x,y,t)\in U^{[1]}$
such that $t\not=0$.}
\end{equation}
If, conversely, $f\!: U\to F$ is continuous
and (\ref{deff1}) holds for
a continuous map $f^{[1]}\!: U^{[1]}\to F$,
then $f$ is $C^1$, with $df(x,y)=\lim_{t\to 0}t^{-1}(f(x+ty)-f(x))=\lim_{t\to0}
f^{[1]}(x,y,t)=f^{[1]}(x,y,0)$.
\end{numba}
The preceding characterization of $C^1$-maps
is a useful tool for various purposes.
Beyond the real case,
the characterizing property
just described can be used to {\em define\/} $C^1$-maps~\cite{BGN}:
\begin{numba}\label{defnBGN}
Let $E$ and $F$ be (Hausdorff) topological vector
spaces over a topological field
$\K$ (which we always assume
Hausdorff and non-discrete),
and $U\sub E$ be open. Let $U^{[1]}:=\{(x,y,t)\in U\times E\times\K\!:
x+ty\in U\}$. A map $f\!: U\to F$ is called $C^1$
if it is continuous and there exists a (necessarily unique)
continuous map $f^{[1]}\!: U^{[1]}\to F$
such that (\ref{deff1}) holds.
Inductively, $f$ is called $C^{k+1}$
for $k\in \N$ if $f$ is $C^1$ and $f^{[1]}\!: U^{[1]}\to F$
is $C^k$. The map $f$ is $C^\infty$ or {\em smooth\/}
if it is $C^k$ for all $k\in \N$. We write
$C^k_\K$ for $C^k$ if we wish to emphasize the ground field.
\end{numba}
By \cite[Prop.\,7.4]{BGN}, the definitions of $C^k$-maps
given in {\bf \ref{defnMB}} and {\bf \ref{defnBGN}} are equivalent
for maps into real locally convex spaces.
Compositions of $C^k$-maps being~$C^k$
\cite[Prop.\,4.5]{BGN},
manifolds and (smooth) Lie groups
modelled on topological $\K$-vector spaces
can be defined in the usual way.
For further information,
see \cite{Mil} (real case)
and~\cite{BGN}. Examples of infinite-dimensional
Lie groups over topological fields can be found in~\cite{ZOO}.
\begin{numba}
A {\em valued field\/}
is a field~$\K$, equipped with an absolute value
$|.|\!: \K\to [0,\infty[$ (see \cite{Wie});
we require furthermore that the absolute value be non-trivial
(i.e., the corresponding metric defines
a non-discrete topology on~$\K$).
Every valued field is, in particular,
a topological field.
A topological vector space~$E$ over a valued field~$\K$
is called {\em polynormed\/}
if its vector topology arises from a family
of continuous seminorms $q\!: E\to [0,\infty[$.
Thus polynormed vector spaces over $\K\in \{\R,\C\}$
are the usual locally convex spaces.
We also write $\|.\|_q:=q$, for better readability.
Given $x\in E$ and $r>0$, we let
$B_r^q(x):=\{y\in E\!: \|y-x\|_q<r\}$
be the open $q$-ball of radius~$r$ around $x$.
\end{numba}
Our studies hinge
on Taylor's formula \cite[Thm.\,5.1]{BGN}:
\begin{prop}\label{basictaylor}
If $k\in \N$ and $f\!: E\supseteq U\to F$ is $C^k$,
then there are
continuous functions
$a_j\!: U\times E\to F$ for $j=1,\ldots, k$
and a continuous function $R_k\!: U^{[1]}\to F$
such that
\[
f(x+ty)-f(x)=\sum_{j=1}^k t^j a_j(x,y)\;+\, t^k\,R_k(x,y,t)\qquad
\mbox{for all $(x,y,t)\in U^{[1]}$}
\]
and $R_k(x,y,0)=0$ for all $(x,y)\in U\times E$.
The functions $a_j$ and $R_k$ are uniquely determined,
$a_j(x,\sbull)$ is homogeneous of degree~$j$,
and $j! a_j(x,y)= d^jf(x,y,\ldots,y)$
for all $(x,y)\in U\times E$.\Punkt
\end{prop}
Here $d^jf\!: U\times E^j\to F$ denotes the $j$th differential of $f$,
defined in terms of
iterated directional derivatives via
$d^jf(x,y_1,\ldots,y_j):=(D_{y_1}\cdots D_{y_j}f)(x)$.
\begin{la}\label{secorder}
Let $E$ and $F$ be polynormed
vector spaces over a valued field~$\K$
and
$f\!: U\to F$ be a $C^2$-map on an open subset $U\sub E$.
Let $x_0\in U$,
$q$ be a continuous seminorm
on~$F$, and $C>0$.
Then there exists a continuous seminorm $p$
on~$E$ such that
$B^p_2(x_0)\sub U$ and
$\|f(x+y)-f(x)-df(x,y)\|_q=
\|R_1(x,y,1)\|_q\leq C\, \|y\|_p^2$ for all $x\in B^p_1(x_0)$
and $y\in B^p_1(0)$.\vspace{1.5mm}
\end{la}
\begin{center}
{\bf H\"{o}lder continuity at a point}\vspace{-1mm}
\end{center}
Until {\bf \ref{endvalued}},
$\K$ denotes a valued field.
\begin{defn}
Let $E$ and $F$ be polynormed
$\K$-vector spaces, $x\in E$,
$U\sub E$ be a neighbourhood of~$x$,
$f\!: U\to F$ be a map,
and $\alpha\in \;]0,1]$.
We say that $f$ is {\em H\"{o}lder continuous
of degree\/} (or H\"{o}lder exponent) {\em $\alpha$
at~$x$\/}
(for short: $f$ is $H_\alpha$ at $x$)
if, for every continuous seminorm $q$ on $F$,
there exist $\delta>0$, $C>0$ and a continuous
seminorm $p$ on $E$ such that $B^p_\delta(x)\sub U$
and
\begin{equation}\label{hcondit}
\|f(y)-f(x)\|_q\leq C\, (\|y-x\|_p)^\alpha\qquad
\mbox{for all $y\in B^p_\delta(x)$.}
\end{equation}
If $f$ is $H_1$ at~$x$, we also say that $f$ is
{\em Lipschitz continuous at~$x$.}
We say that $f$ is {\em H\"{o}lder continuous at~$x$\/}
if $f$ is $H_\alpha$ at~$x$ for some $\alpha\in\,]0,1]$.
\end{defn}
\begin{rem}
Replacing $p$ with
$\max \big\{\delta^{-1}, C^{\frac{1}{\alpha}}\big\}\cdot p$,
we can always achieve that $C=\delta=1$.
\end{rem}

\begin{la}\label{basehold}
For maps between subsets of polynormed
$\K$-vector spaces, we have:
\begin{itemize}
\item[\rm (a)]
If $f$ is $H_\alpha$ at $x$ then $f$ is continuous at $x$.
\item[\rm (b)]
If $\alpha\geq \beta$ and $f$ is $H_\alpha$ at $x$, then
$f$ is $H_\beta$ at $x$.
\item[\rm (c)]
Any $C^1$-map is Lipschitz continuous at each
point.
\item[\rm (d)]
If $f$ is $H_\alpha$ at $x$ and $g$ is $H_\beta$ at $f(x)$,
then $g\circ f$ is $H_{\alpha\cdot \beta}$ at $x$. 
\end{itemize}
\end{la}
\begin{defn}
Let $f\!: M\to N$ be a map between $C^1_\K$-manifolds
modelled on polynormed $\K$-vector spaces, and $\alpha\in \;]0,1]$.
We say that $f$ is {\em H\"{o}lder continuous of degree
$\alpha$ at $x\in M$\/}
(or briefly: $F$ is $H_\alpha$ at~$x$),
if $f$ is continuous at $x$ and there are
a chart
$\phi\!: U_1\to U$ of $M$ around~$x$
and a chart $\psi\!: V_1\to V$ of $N$ around
$f(x)$, such that
$\,\phi(f^{-1}(V_1)\cap U_1)\to
V$,
$y\mto \psi(f(\phi^{-1}(y)))\,$
is $H_\alpha$ at $\phi(x)$.
\,(This then holds for any choice of $\phi$ and $\psi$, by
La.\,\ref{basehold}).\vspace{1mm}
\end{defn}
\begin{center}
{\bf Notions of differentiability at a point}\vspace{-1mm}
\end{center}
\begin{numba}\label{dftot}
(Cf.\ \cite[I,\,\S3]{Lan}).
Let $E$ and $F$ be topological $\K$-vector spaces,
$x\in E$, and $f\!: U\to F$ be a map
defined on a neighbourhood~$U$ of~$x$ in~$E$.
The map $f$ is called {\em totally differentiable at $x$\/} if
there is a  (necessarily unique) continuous linear map
$f'(x)\!: E\to F$ such that
\[
h\!: U-x\to F,\qquad h(y):= f(x+y)-f(x)-f'(x).y
\]
is {\em tangent to~$0$\/}
in the sense that,
for every $0$-neighbourhood $W\sub F$,
there is a $0$-neighbourhood $V\sub E$
and a function
$\theta\!: I\to \K$ defined on some
$0$-neighbourhood~$I\sub \K$
such that $I\cdot V\sub U-x$,
$\theta(t)=o(t)$ (i.e.,
$\theta(0)=0$ and
$\lim_{t\to 0}\theta(t)/t=0$),
and
\[
h(tV)\sub \theta(t)W \qquad\mbox{for all $\,t\in I$.}
\]
\end{numba}
\begin{numba}\label{claim1tot}
If $E$ and $F$ are polynormed, then
$h$ as before
is tangent to $0$ if and only if,
for every continuous seminorm $q$ on~$F$,
there exists a continuous seminorm $p$ on $E$
such that, for each $\ve>0$, there exists
$\delta>0$ such that $B_\delta^p(0)\sub U-x$ and
\[
\|h(y)\|_q\leq \ve\|y\|_p\qquad\mbox{for all $y\in B^p_\delta(0)$.}
\]
\end{numba}
\begin{numba}\label{claim2tot}
The Chain Rule holds:
If $f\!: E\supseteq U\to F$
is totally differentiable at~$x$
and the map $g\!: F\supseteq V\to H$ is totally differentiable
at $f(x)$ and $f(U)\sub V$, then $g\circ f\!: U\to H$
is totally differentiable at~$x$, with $(g\circ f)'(x)=g'(f(x))\circ
f'(x)$.
\end{numba}
\begin{la}\label{C2total}
Let $E$ and $F$ be topological
$\K$-vector spaces, $U\sub E$ be an open subset, and
$f\!: U\to F$ be a $C^2$-map.
Then $f$ is totally differentiable at
each $x\in U$,
and $f'(x)=df(x,\sbull)$.
\end{la}
\begin{numba}\label{endvalued}
Given $r\in \N\cup\{\infty\}$,
a map $f\!: M\to N$
between $C^r$-manifolds modelled on topological
$\K$-vector spaces,
and $x\in M$, we call $f$ {\em totally differentiable
at $x$\/} if $f$ is continuous at $x$
and there exist a chart
$\phi\!: U_1\to U$ of $M$ around~$x$
and a chart $\psi\!: V_1\to V$ of $N$ around
$f(x)$, such that
$\,\phi(f^{-1}(V_1)\cap U_1)\to
V$,
$y\mto \psi(f(\phi^{-1}(y)))\,$
is totally differentiable at
$\phi(x)$.\footnote{If $r\geq 2$,
then the latter property
is independent of the choice of charts,
by the Chain Rule
(the chart changes
are $C^2$ and hence totally differentiable at each point
by Lemma~\ref{C2total}).}
\end{numba}
We find it convenient to work with
a certain weaker differentiability property,
which even makes sense over
arbitrary topological fields:
\begin{numba}\label{defnfeeb}
Let $E$ and $F$ be topological vector spaces
over a topological field~$\K$,
$U\sub E$ be open, $x\in U$, and
$f\!: U\to F$ a continuous map.
Let $A:=\{(y,t)\in E\times \K^\times\!: x+ty\in U\}$
and $\wt{U}_x:=A \cup (E\times \{0\})\sub E\times \K$.
We say that $f$ is {\em feebly differentiable
at~$x$\/} if there is a (unique)
continuous linear
map $f'(x)\!: E\to F$
making the following map continuous:
\[
\wt{f}_x\!: \wt{U}_{x}\to F,\qquad
(y,t)\mto \left\{
\begin{array}{cl}
\frac{f(x+ty)-f(x)}{t} &\mbox{if $\;t\not=0$}\\
f'(x).y &\mbox{if $\;t=0$.}
\end{array}
\right.
\]
\end{numba}
\begin{la}\label{2variable}
Let $E$ and $F$ be topological vector spaces
over a topological field~$\K$,
$U\sub E$ be open,
$f\!: U\to F$ be a map, and
$x\in U$.
If $f$ is $C^1$ or if
$\K$ is a valued field,
$f$ is continuous on~$U$
and totally differentiable
at~$x$, then $f$ is feebly differentiable at~$x$.
\end{la}
\begin{numba}\label{chainfeeb}
The Chain Rule holds for feebly differentiable
maps: If $f\!: E\supseteq U\to F$
is feebly differentiable at~$x$
and $g\!: F\supseteq V\to H$ is feebly differentiable
at $f(x)$ and $f(U)\sub V$, then $g\circ f\!: U\to H$
is feebly differentiable at~$x$, with $(g\circ f)'(x)=g'(f(x))\circ
f'(x)$.
\end{numba}
\begin{numba}
A map $f\!: M\to N$ between $C^1$-manifolds
modelled on topological $\K$-vector spaces is called
{\em feebly differentiable at $x\in M$\/}
if it is continuous at~$x$
and $y\mto \psi(f(\phi^{-1}(y)))$
is feebly differentiable at
$\phi(x)$
for charts $\phi$ and $\psi$ as in {\bf \ref{endvalued}}.
\end{numba}
Cf.\ \cite{AaS} for a comparative study of various
differentiability properties
at a point.
\section{Homomorphisms between Lie groups}\label{sechom}
We prove preparatory results
concerning differentiability properties of homomorphisms.
\begin{la}\label{C1smooth}
Let $f\!: G\to H$ be a $C^1_\K$-homomorphism
between
Lie groups over $\K\in \{\R,\C\}$,
where $H$ is modelled on a locally convex space.
Then $f$ is $C^\infty_\K$.
\end{la}
\begin{proof}
We show that $f$ is
$C^k$ for each $k\in \N$, by induction.
By hypothesis, $f$ is $C^1$.
Using the
trivialization
$\tau_G\!: G\times L(G)\to TG$,
$\tau_G(g,X):=T_1\lambda_g(X)$
(where $\lambda_g\!: G\to G$, $x\mto gx$
denotes left translation by~$g$)
and the corresponding trivialization
$\tau_H\!: H\times L(H)\to TH$,
the tangent map $Tf$
can be expressed as
\begin{equation}\label{givesindu}
Tf=\tau_H\circ (f\times L(f))\circ (\tau_G)^{-1}\,.
\end{equation}
Since $\tau_G$ and $\tau_H$ are $C^\infty$-diffeomorphisms
and the continuous linear map
$L(f)$ is smooth,
(\ref{givesindu}) shows that
if $f$ is $C^k$,
then so is $Tf$.
But then $f$ being a $C^1$-map into a manifold
modelled on a {\em locally convex\/} space
with $Tf$
of class $C^k$, the map $f$ is $C^{k+1}$
(cf.\ \cite[Prop.\,7.4]{BGN}).
\end{proof}
\begin{la}\label{totsmooth}
Let $f\!: G\to H$ be a
homomorphism between
Lie groups modelled on
topological vector spaces over a topological field~$\K$.
Assume that $f$ is feebly differentiable
at~$1$ $($this is the case if $\K$
is a valued field and $f$ is totally differentiable at~$1)$.
Then $f$ is $C^1_\K$.
If
$\K\in \{\R,\C\}$
and the modelling space of $H$ is locally convex,
then $f$ is $C^\infty_\K$.
\end{la}
\begin{proof}
We let $\phi\!: U_1\to U\sub L(H)$
be a chart of~$H$ around~$1$, such that $\phi(1)=0$.
There exists an open
identity neighbourhood $V_1\sub U_1$
such that $V_1V_1\sub U_1$;
let $V:=\phi(V_1)$. Then
\[
\mu\!: V\times V\to U,\qquad \mu(x,y):=x*y:=\phi(\phi^{-1}(x)\phi^{-1}(y))
\]
expresses multiplication on~$H$ in local coordinates.
Let $\psi\!: P_1\to P\sub L(G)$ be a chart of $G$
such that $f(P_1)\sub U_1$, $0\in P$ and $\psi(1)=0$;
let $Q_1$ and $S_1$ be open
identity neighbourhoods in~$G$ such that $Q_1Q_1\sub P_1$,
$f(Q_1)\sub V_1$, $S_1=(S_1)^{-1}$, and $S_1S_1\sub Q_1$.
Then $Q:=\psi(Q_1)$
and $S:=\psi(S_1)$ are open $0$-neighbourhoods
in $L(G)$.
Define $\iota\!: S\to S$,
$\iota(x):=x^{-1}:=\psi(\psi^{-1}(x)^{-1})$
and
$\nu\!: Q\times Q\to P$, $\nu(x,y):=x*y:=\psi(\psi^{-1}(x)\psi^{-1}(y))$.
Then
\[
g:=\phi\circ f|_{P_1}^{U_1}\circ \psi^{-1}\!: P \to U
\]
maps $0$ to $0$
and is continuous (since $f$ is continuous, being a
homomorphism which is continuous at one point).
Furthermore,
$g$ is feebly
differentiable at~$0$
by hypothesis (resp., Lemma~\ref{2variable}).
For $(x,y,t)\in S^{]1[}:=\{(x,y,t)\in S^{[1]}\!: t\not=0\}$,
we have
\begin{eqnarray*}
t^{-1}(g(x+ty)-g(x)) &=& t^{-1}(g(x)*g(x^{-1}*(x+ty))-g(x))\\
&=&
t^{-1}(g(x)*(0+tt^{-1}g(x^{-1}*(x+ty)))-g(x)*0)\\
&=& \mu^{[1]}((g(x),0),\, (0,t^{-1}g(x^{-1}*(x+ty))),\, t)\\
&=& \mu^{[1]}((g(x),0),\, (0,t^{-1}g(th(x,y,t))),\, t)\\
&=& \mu^{[1]}((g(x),0),\, (0,\wt{g}_0(h(x,y,t),t)),\, t)
\end{eqnarray*}
where
$h\!: S^{[1]}\to L(G)$,
$h(x,y,t):=\nu^{[1]}((x^{-1},x),\, (0,y),\, t)$
is continuous, and so is
the map
$\wt{g}_0\!: \wt{P}_0\to L(H)$
(defined as in {\bf \ref{defnfeeb}}).
Note that $F\!: S^{[1]}\to L(H)$,
$F(x,y,t):=\mu^{[1]}((g(x),0),\, (0,\wt{g}_0(h(x,y,t),t)),\, t)$
makes sense on all of $S^{[1]}$.
The map $F$ is continuous and,
by the preceding,
we have $F(x,y,t)=\frac{1}{t}(g(x+ty)-g(x))$
for all $(x,y,t)\in S^{]1[}$.
Thus $g|_S$ is $C^1$,
with $(g|_S)^{[1]}=F$.
Hence $f|_{S_1}$ is $C^1$
and hence so is $f$ on all of $G$,
by \cite[La.\,3.1]{ANA}.
If $\K\in\{\R,\C\}$
and $L(H)$ is locally convex,
this entails that $f$ is $C^\infty$
(Lemma~\ref{C1smooth}).
\end{proof}
\section{H\"{o}lder continuous homomorphisms are smooth}\label{sechold}
In this section, which is the core of the article,
we establish the main result.
\begin{defn}
A sequence $(x_n)_{n\in \N}$
in a topological vector space $E$ over a valued
field $\K$ is called {\em Mackey-Cauchy\/}
if there exists a bounded subset $B\sub E$
and elements $\mu_{n,m}\in \K$ such that
$x_n-x_m\in \mu_{n,m}B$
for all $n,m\in \N$ and $\mu_{n,m}\to 0$ as both
$n,m\to\infty$
(cf.\ \cite[p.\,14]{KaM}).
We say that $E$ is {\em Mackey complete\/}
if every Mackey-Cauchy sequence in~$E$
is convergent (cf.\ \cite[La.\,2.2]{KaM}).
\end{defn}
\begin{thm}\label{main}
Let $f\!:G\to H$ be a homomorphism between
smooth Lie groups modelled on
locally convex, real topological vector spaces.
If
the modelling space of $H$ is Mackey complete
and $f$ is H\"{o}lder continuous at~$1$,
then $f$ is smooth.
\end{thm}
\begin{proof}
By hypothesis, $f$ is $H_\alpha$
at~$1$ for some $\alpha\in \;]0,1]$.
The proof proceeds in two steps.
The first goal is to show that
if $\alpha\in \;]0,\half]$,
then $f$ also is
$H_{\frac{3}{2}\alpha}$ at~$1$.
Since the H\"{o}lder exponent
can be improved repeatedly, this means that $f$ actually
is
$H_\alpha$ at~$1$ with $\alpha\in \;]\half,1]$.
Having achieved this,
the second goal will be
to show that $f$ is totally differentiable
at~$1$ and hence smooth, by Lemma~\ref{totsmooth}.\\[3mm]
For the moment, we only know that $\alpha\in\;]0,1]$.
We let $\phi\!: U_1\to U\sub L(H)$
be a chart of~$H$ around~$1$, such that $\phi(1)=0$.
There exist open, symmetric\footnote{Recall that an identity neighbourhood
$X$ is {\em symmetric\/} if $X=X^{-1}$.} 
identity neighbourhoods $V_1\sub U_1$ and $W_1\sub V_1$
such that $V_1V_1\sub U_1$ and $W_1W_1\sub V_1$;
let $V:=\phi(V_1)$ and $W:=\phi(W_1)$. Then
\[
\mu\!: V\times V\to U,\qquad \mu(x,y):=x*y:=\phi(\phi^{-1}(x)\phi^{-1}(y))
\]
expresses the multiplication of~$H$ in local coordinates.
Products of more than two elements are formed from left to right;
for example, $x*y*z:=(x*y)*z$. Of course, $(x*y)*z=x*(y*z)$
whenever both products are defined (and likewise
for products of more than three factors). 
Since $0*0=0$ and $\mu'(0,0).(x,y)=x+y$, the map
\[
\sigma\!: W\times W \to U,\quad \sigma(x,y):=x*x*y
\]
satisfies $\sigma(0,0)=0$ and $\sigma'(0,0)(u,v)=2u+v$
for $u,v\in L(H)$.
Hence, using the Taylor expansion of
$\sigma$ about $(0,0)$, we have
\[
\sigma(x,y)=2x+y+R(x,y)\qquad\mbox{for all $x,y\in W$,}
\]
where $R(x,y):=R_1((0,0),(x,y),1)$ (cf.\ Proposition~\ref{basictaylor}).
Let $\psi\!: P_1\to P\sub L(G)$ be a chart of $G$
around~$1$, such that $f(P_1)\sub U_1$ and $\psi(1)=0$;
let $Q_1\sub P_1$ and $B_1\sub Q_1$ be symmetric
identity neighbourhoods such that $Q_1Q_1\sub P_1$, $B_1B_1\sub Q_1$,
$f(Q_1)\sub V_1$, and $f(B_1)\sub W_1$.
Set $Q:=\psi(Q_1)$ and $B:=\psi(B_1)$.
Define $\iota\!: Q\to Q$,
$\iota(x):=x^{-1}:=\psi(\psi^{-1}(x)^{-1})$
and
$\nu\!: Q\times Q\to P$, $\nu(x,y):=x*y:=\psi(\psi^{-1}(x)\psi^{-1}(y))$.
Then
\[
g:=\phi\circ f|_{P_1}^{U_1}\circ \psi^{-1}\!: P \to U
\]
is continuous, maps $0$ to $0$,
and is $H_\alpha$ at~$0$.\\[3mm]
We now adapt the ideas explained in the Introduction
for the special case of one-parameter groups
to the present, fully general situation.
To this end, let
$A\sub B$ be a balanced, open
$0$-neighbourhood such that $A*A\sub
B$,
$\iota(A)*\iota(A)\sub B$,
and $\iota(A)*\iota(A)*A\sub g^{-1}(W)$.
We abbreviate
$(\half x)^{-2}:=\iota(\half x)*\iota(\half x)$
for $x\in A$ and define
\begin{equation}\label{defnh}
h\!: A\to W,\quad h(x):=
g((\half x)^{-2}* x)\, .
\end{equation}
We have
$g(x)= g((\half x)^2*(\half x)^{-2}*x)=
g(\half x)^2* g((\half x)^{-2}* x)
=\sigma(g(\half x),g((\half x)^{-2}x))
=2 g(\half x) + g((\half x)^{-2}x)
+R(g(\half x),g((\half x)^{-2}x))$
for $x\in A$ and hence
\begin{equation}\label{start}
g(\half x)=\half g(x)
-\half h(x)
-\half
R(g(\half x),h(x))\,,
\end{equation}
with $h$ as in (\ref{defnh}).
Since also $\half x\in A$,
likewise
$g(\quart x) = \half g(\half x)
-\half h(\half x)
-\half
R(g(\quart x),h(\half x))$.
Inserting the right hand side of (\ref{start}) for
$g(\half x)$ here, we arrive at
\[
g(\quart x)=
\quart g(x)
-\quart h(x)
-\quart
R(g(\half x),h(x)) -\half h(\half x)
-\half
R(g(\quart x),h(\half x))\,.
\]
Proceeding in this way, we obtain
\begin{equation}\label{trick1}
g(2^{-n}x)=
2^{-n}g(x)-\sum_{k=1}^n
2^{-n+k-1}
\left(
h(2^{1-k}x)+
R(g(2^{-k}x),
h(2^{1-k}x))\right)
\end{equation}
for all $n\in \N_0$, by induction.
Hence
\begin{equation}\label{trick2}
2^n g(2^{-n}x)=
g(x)-\sum_{k=1}^n
\,2^{k-1}
\left(
h(2^{1-k}x)+
R(g(2^{-k}x),
h(2^{1-k}x))\right)
\end{equation}
for all $x\in A$ and
$n\in \N_0$.
The following lemma provides
estimates on the summands in (\ref{trick2});
later, these estimates will be used
to show that the series is summable
(see (\ref{thesum})).
\begin{la}\label{LaNull}
Let $q$ be a continuous seminorm on $L(H)$.
Then there exists a continuous seminorm
$p$ on $L(G)$
such that $B_1^p(0) \sub A$,
\begin{equation}\label{desired}
\|h(x)+R(g(\half x), h(x))\|_q\leq \|x\|_p^{2\alpha}\qquad
\mbox{for all $\, x\in B_1^p(0)$,}
\end{equation}
and $\|g(x)\|_q\leq (\|x\|_p)^\alpha$ for all $x\in B_1^p(0)$.
\end{la}
\begin{proof}
As a consequence of Lemma~\ref{secorder},
there exists a continuous
seminorm $r$ on $L(H)$
such that $B_1^r(0)\sub W$ and
\begin{equation}\label{usetaylor}
\|R(y,z)\|_q\leq \half (\max\{\|y\|_r,\|z\|_r\})^2
\quad
\mbox{for all $y,z\in B_1^r(0)\,$;}
\end{equation}
after replacing $r$ with $r+q$,
we may assume that $r\geq q$.
Since $g$ is
$H_\alpha$ at~$0$,
there is a continuous seminorm
$s$ on $L(G)$
such that $B_1^s(0)\sub P$,
$g(B_1^s(0))\sub B_1^r(0)$,
and
\begin{equation}\label{usehold}
\|g(x)\|_r\leq \half (\|x\|_s)^\alpha\qquad \mbox{for all $x\in B_1^s(0)$.}
\end{equation}
We now consider the smooth
map $j\!: A\to Q$, $j(x):=(\frac{1}{2}x)^{-2}* x$.
Then $j(0)=0$ and $j'(0)=0$,
entailing that there exists
a continuous seminorm $p$
on $L(G)$ such that $B_1^p(0)\sub A$,
$j(B_1^p(0))\sub B_1^s(0)$,
and
\begin{equation}\label{appltay}
\|j(x)\|_s\leq (\|x\|_p)^2\qquad \mbox{for all $x\in B_1^p(0)$}
\end{equation}
(cf.\ Lemma~\ref{secorder});
we may assume that $p\geq s$.
Then
$\|h(x)\|_q\leq \|h(x)\|_r=
\|g(j(x))\|_r\leq \half (\|j(x)\|_s)^\alpha
\leq \half (\|x\|_p)^{2\alpha}$
for all $x\in B_1^p(0)$,
by (\ref{usehold}) and (\ref{appltay}).
Also
$\|g(\half x)\|_r\leq \half (\|\half x\|_s)^\alpha
\leq (\|x\|_s)^\alpha\leq (\|x\|_p)^\alpha$
and $\|h(x)\|_r\leq \half (\|x\|_p)^{2\alpha} \leq (\|x\|_p)^\alpha$,
whence
$\|R(g(\half x),h(x))\|_q\leq \half (\|x\|_p)^{2\alpha}$,
by (\ref{usetaylor}).
Using the preceding estimates, we obtain
$\|h(x)+R(g(\half x),h(x))\|_q\leq
\|h(x)\|_q+\|R(g(\half x),h(x))\|_q\leq
\half (\|x\|_p)^{2\alpha}+ \half (\|x\|_p)^{2\alpha}
=(\|x\|_p)^{2\alpha}$ for all
$x\in B_1^p(0)$.
Thus (\ref{desired}) holds.
We also have
$\|g(x)\|_q\leq\|g(x)\|_r\leq \half (\|x\|_s)^\alpha
\leq (\|x\|_s)^\alpha\leq (\|x\|_p)^\alpha$.
\end{proof}
\begin{la}\label{impalph}
If $f$ is $H_\alpha$
at $1$ with
$\alpha\in \;]0,\frac{1}{2}]$,
then
$f$ also is
$H_{\frac{3}{2}\alpha}$ at~$1$.
\end{la}
\begin{proof}
Given a continuous seminorm $q$ on $L(H)$,
we let $p$ be as in Lemma~\ref{LaNull}.
In the following, we show that
\begin{equation}\label{givhol}
\|g(y)\|_q \leq K\,
2^{\frac{3}{2}\alpha}(\|y\|_p)^{\frac{3}{2}\alpha}
\quad\mbox{for all $y\in B^p_1(0)$,}
\end{equation}
for a suitable constant $K\in [0, \infty[$.
Thus $g$ will be
$H_{\frac{3}{2}\alpha}$ at~$0$, and hence
$f$ will be
$H_{\frac{3}{2}\alpha}$ at~$1$.

Using (\ref{trick1}) and the estimates from
Lemma~\ref{LaNull}, we obtain
\begin{eqnarray}
\|g(2^{-n}x)\|_q &\leq &
2^{-n}\|g(x)\|_q+\sum_{k=1}^n 2^{-n+k-1}\|h(2^{1-k}x)
+R(g(2^{-k}x),h(2^{1-k}x))\|_q\nonumber\\
&\leq &
2^{-n} +\sum_{k=1}^n
2^{-n+k-1}(\|2^{1-k}x\|_p)^{2\alpha}\leq
2^{-n} + \sum_{k=1}^n 2^{-n+k-1}
2^{2\alpha-2\alpha k}\nonumber \\
&=&
\left(2^{-(1-\frac{3}{2}\alpha)n}
+  2^{2\alpha-1}
2^{-(1-\frac{3}{2}\alpha)n}
\sum_{k=1}^n
2^{(1-2\alpha)k}\right) 2^{-\frac{3}{2}\alpha n}\label{getimpro}
\end{eqnarray}
for all $x\in B_1^p(0)$ and $n\in \N_0$.
Here $2^{-(1-\frac{3}{2}\alpha)n}\leq 1$
for all $n\in \N_0$.
The summation formula
for the finite geometric series yields
\[
\sum_{k=1}^n2^{(1-2\alpha)k}
\leq
\sum_{k=1}^n2^{(1-\frac{3}{2}\alpha)k}=
\frac{2^{(1-\frac32 \alpha)(n+1)}-2^{1-\frac32 \alpha}}{2^{1-\frac32 \alpha}-1}
\leq
\frac{2^{(1-\frac32 \alpha)(n+1)}}{2^{1-\frac32 \alpha}-1}
=c\, 2^{(1-\frac32\alpha) n}
\]
with $c:=\frac{2^{1-\frac32 \alpha}}{2^{1-\frac32 \alpha}-1}$.
We therefore obtain the following estimates
for the second term
in (\ref{getimpro}):
\[
2^{2\alpha-1}
2^{-(1-\frac{3}{2}\alpha)n}
\sum_{k=1}^n
2^{(1-2\alpha)k}
\leq
c \,
2^{2\alpha-1}
2^{-(1-\frac{3}{2}\alpha)n}2^{(1-\frac32 \alpha) n}
=K_1
\]
for all $n\in \N_0$,
with $K_1:=
c \,
2^{2\alpha-1}$. Using the estimates
just established, (\ref{getimpro}) yields
\begin{equation}\label{impro2}
\|g(2^{-n}x)\|_q\leq
K\, 2^{-\frac{3}{2}\alpha n}\qquad
\mbox{for all $x\in B^p_1(0)$ and $n\in \N_0$,}
\end{equation}
with $K:=1+K_1$.
Then (\ref{givhol}) holds with $K$ as just defined.
To see this, let $y\in B^p_1(0)$.
If $\|y\|_p=0$, then $\|g(y)\|_q\leq \|y\|_p^\alpha=0
\leq K 2^{\frac{3}{2}\alpha} \|y\|_p^{\frac{3}{2}\alpha}$,
as desired.
If $\|y\|_p>0$,
then there exists $n\in \N_0$ such that
$2^{-n-1} \leq \|y\|_p<2^{-n}$.
Thus $2^{-n}\leq 2 \|y\|_p$.
Since
$y=2^{-n}x$
with
$x:=2^ny \in B^p_1(0)$,
(\ref{impro2}) yields
\[
\|g(y)\|_q=\|g(2^{-n}x)\|_q
\leq K (2^{-n})^{\frac{3}{2}\alpha}
\leq K(
2 \|y\|_p)^{\frac{3}{2}\alpha}
=
K 2^{\frac{3}{2}\alpha}
(\|y\|_p)^{\frac{3}{2}\alpha}\,,
\]
whence (\ref{givhol})
also holds if $\|y\|_p>0$. This completes the proof of Lemma~\ref{impalph}.
\end{proof}
If $\alpha\in \;]0,\half]$,
there exists $k\in \N$ such that
$(\frac{3}{2})^{k-1}\alpha\leq \half$
and $\beta:=(\frac{3}{2})^k\alpha
\in \;] \half,1]$.
Repeated application of Lemma~\ref{impalph}
shows that
$f$ is
$H_\beta$ at~$1$.
After replacing $\alpha$
with $\beta$,
we may assume
throughout the following that
$\alpha\in \;]\half,1]$.\\[3mm]
In the remainder of the proof,
we show that
$g$ is totally
differentiable at~$0$. The main point
is to construct a candidate $\Lambda$ for
the derivative $g'(0)$.
We first construct $\lambda=\Lambda|_A$
on the $0$-neighbourhood $A\sub L(G)$ (from above).
\begin{la}\label{La1}
The limit
$\lambda(x):=\lim_{n\to\infty} \frac{g(2^{-n}x)}{2^{-n}}$
exists in $L(H)$, for each $x\in A$.
For each continuous seminorm
$q$ on $L(H)$,
the convergence of
$\frac{g(2^{-n}x)}{2^{-n}}$
in $(L(H),\|.\|_q)$ is locally uniform
in~$x$. The map $\lambda\!: A\to L(H)$ is continuous.
\end{la}
\begin{proof}
Fix $x_0\in A$.
Given a continuous seminorm
$q$ on
$L(H)$, we let
$p$ be as in Lemma~\ref{LaNull}.
There is $N\in \N$ such that $2^{-N}\|x_0\|_p<1$.
Then $S:=B^p_{2^N}(0)\cap A$
is an open neighbourhood of $x_0$ in $A$ such that
$2^{-N}S\sub B^p_1(0)\sub A$.
Abbreviate $C:= 2^{2\alpha N}$
and $K:=\frac{C2^{2\alpha-1}}{1-2^{-(2\alpha-1)}}$.
Let $M\geq N$.
For all $m,n\geq M$ (where $m\geq n$, say),
using (\ref{trick2})
we obtain for all $x\in S$:
\begin{eqnarray}
\!\!\left\|2^mg(2^{-m}x)
-2^n g(2^{-n}x)\right\|_q
&=&
\left\| \sum_{k=n+1}^m 2^{k-1}
\left(
h(2^{1-k}x)+
R(g(2^{-k}x),
h(2^{1-k}x))\right)\right\|_q\nonumber\\
&\leq &
\sum_{k=n+1}^m 2^{k-1}
\|h(2^{1-k}x)+
R(g(2^{-k}x),
h(2^{1-k}x))\|_q\nonumber\\
&\leq&
\underbrace{\|x\|_p^{2\alpha}}_{\leq C}
\sum_{k=n+1}^m2^{k-1}2^{2\alpha(1-k)}
\leq
C\, 2^{2\alpha-1}\!\! \sum_{k=n+1}^m 2^{-(2\alpha-1)k}\label{auchbrauch}\\
&\leq & K\cdot (2^{-(2\alpha-1)})^{n+1}
\leq
K\cdot ( 2^{-(2\alpha-1)})^{M+1}\,,\label{unf}
\end{eqnarray}
using (\ref{desired}) to pass to the
third line, then using that $2^{-(2\alpha-1)}<1$
since $\alpha \in \;]\frac{1}{2},1]$.
Here, the final expression tends to $0$ as $M\to\infty$,
uniformly in $x \in S$.

By the preceding considerations,
$\left(2^ng(2^{-n}x_0)\right)_{n\in \N_0}$
is a Cauchy sequence in $L(H)$
in particular. Hence, if $L(H)$ is sequentially
complete, the limit
\begin{equation}\label{thesum}
\lambda(x_0):=\lim_{n\to\infty} 2^n g(2^{-n}x_0)
=g(x_0)-\sum_{k=1}^\infty 2^{k-1}(h(2^{1-k}x_0)+
R(g(2^{-k}x_0),h(2^{1-k}x_0)))
\end{equation}
exists in $L(H)$. As we shall presently see,
the limit also exists when $L(H)$ is Mackey complete.
Assuming the validity of this claim
for the moment,
letting $m\to \infty$ in the lines before (\ref{unf})
we obtain
$\big\|\lambda(x)-2^ng(2^{-n}x)\big\|_q\leq
K \cdot
( 2^{-(2\alpha-1)})^{M+1}$
for all $n\geq M$.
Hence $\big\|\lambda(x)-2^ng(2^{-n}x)\big\|_q\to 0$
uniformly
in $x\in S$, proving the second assertion
of the lemma.
The preceding also entails that $\lambda$ is continuous.\\[3mm]
To complete the proof,
it only remains to prove our claim that
the limit (\ref{thesum}) exists.
Since $L(H)$ is Mackey complete,
we only need to show that $(v_n)_{n\in \N}$
is a Mackey-Cauchy sequence,
where
$v_n:=2^ng(2^{-n}x_0)$.
To this end, pick $a\in \;]2^{-(2\alpha-1)},1[$
and define $r_{n,m}:=a^{\min\{n,m\}+1}$.
Then $r_{n,m}\to 0$ as both $n,m\to\infty$,
and
\[
v_n-v_m\in r_{n,m}\, \Omega \qquad\mbox{for all $n,m\in \N\,$,}
\]
where $\Omega:=\{ r_{n,m}^{-1}(v_n-v_m)\!: n,m\in \N\}$.
If we can show that $\Omega$ is bounded in~$E$,
then $(v_n)_{n\in \N}$ will be Mackey-Cauchy.
To prove boundedness,
assume that
$q$ is a continuous seminorm on $L(H)$.
Let $p$, $N$ and $K$ be as before.
For all $n,m\in \N$,
we have, abbreviating
$\ell:=\max\{N+1,\min\{n,m\}+1\}$:
\begin{eqnarray*}
\|r_{n,m}^{-1}(v_n-v_m)\|_q
&\leq &
a^{-\min\{n,m\}-1}\!\!\!\sum_{k=\min\{n,m\}+1}^{\max\{n,m\}}
\!\!\!\!2^{k-1}  \|h(2^{1-k}x_0)+
R(g(2^{-k}x_0), h(2^{1-k}x_0))\|_q\\
&\leq &
C_q \, +\, a^{-\min\{n,m\}-1}\!\!\!
\sum_{k=\ell}^{\max\{n,m\}}
\!\! \! 2^{k-1} \|h(2^{1-k}x_0)+
R(g(2^{-k}x_0), h(2^{1-k}x_0))\|_q\\
&\leq &
C_q\, +\, a^{-\min\{n,m\}-1}K ( 2^{-(2\alpha-1)})^\ell\\
&\leq & C_q \, +\,  K(a^{-1}2^{-(2\alpha-1)})^{\min\{n,m\}+1}
\leq C_q+K\,,
\end{eqnarray*}
where $C_q:=a^{-N-1}\sum_{k=2}^N 2^{k-1}\|h(2^{1-k}x_0)+
R(g(2^{-k}x_0),
h(2^{1-k}x_0))\|_q$
is an upper bound for
the sum of all terms with $k\leq N$,
for which we do not have estimates available.
Passing to the third line,
we tackled the
summands with $k>N$
as in the proof of (\ref{unf}).
The final inequality holds because
$a^{-1}2^{-(2\alpha-1)}<1$, by the choice of~$a$.
Thus $\|v\|_q\leq C_q+K$ for all $v\in \Omega$,
entailing that $\Omega$
is indeed bounded.
\end{proof}
Before we can prove that $\lambda$
extends to a continuous linear map,
we need another technical result analogous to
Lemma~\ref{LaNull}.\\[3mm]
Let $Z\sub A$ be an open $0$-neighbourhood
such that $Z+Z\sub A$.
We define
$j\!: Z\times Z\to Q$,
$j(x,y):=y^{-1}*x^{-1}*(x+y)$.
Then $j(Z\times Z)\sub g^{-1}(W)$.
The map $\tau\!:
W\times W\times W\to U$, $\tau(x,y,z):=x*y*z$
is smooth, with $\tau(0,0,0)=0$ and
$\tau'(0,0,0)(u,v,w)=u+v+w$
for all $u,v,w\in L(H)$.
Let
$\wt{R}_1\!: (W\times W\times W)^{[1]}\to L(H)$
be the
first order Taylor remainder
of $\tau$.
Abbreviating $D(x,y,z):=\wt{R}_1((0,0,0),\, (x,y,z),\, 1)$,
we then have
\begin{equation}\label{taytau}
\tau(x,y,z)=x+y+z + D(x,y,z)\qquad\mbox{for all $x,y,z\in W$.}
\end{equation}
\begin{la}\label{Latech}
For every continuous seminorm $q$
on $L(H)$,
there is a continuous seminorm~$p$
on $L(G)$
such that $B_1^p(0) \sub Z$ and
\[
\left\|\,  g(j(x,y))\, + \, D \bigl( g(x),\,
g(y),\, g(j(x,y)) \bigr) \, \right\|_q \leq (\max\{\|x\|_p,
\|y\|_p\})^{2\alpha}\quad
\mbox{{\em for all $x,y\in B_1^p(0)$.}}
\]
\end{la}
\begin{proof}
There exists a continuous
seminorm $r$ on $L(H)$
such that
$B_1^r(0)\sub W$,
\begin{equation}\label{usetaylr}
\|D(x,y,z)\|_q\leq \half (\max\{\|x\|_r,\|y\|_r,\|z\|_r\})^2
\quad
\mbox{for all $x,y,z\in B_1^r(0)\,$,}
\end{equation}
and $r\geq q$ (cf.\ Lemma~\ref{secorder}).
Since $g$ is
$H_\alpha$ at~$0$,
there exists a continuous seminorm
$s$ on $L(G)$
such that $B_1^s(0)\sub P$,
$g(B_1^s(0))\sub B_1^r(0)$,
and
\begin{equation}\label{usehld}
\|g(x)\|_r\leq \half (\|x\|_s)^\alpha\qquad
\mbox{for all $x\in B_1^s(0)$.}
\end{equation}
Since $j$ is smooth,
$j(0,0)=0$ and $j'(0,0)=0$,
there exists
a continuous seminorm $p$
on $L(G)$ such that $B_1^p(0)\sub Z$,
$j(B_1^p(0)\times B_1^p(0))\sub B_1^s(0)$,
and
\begin{equation}\label{applty}
\|j(x,y)\|_s\leq (\max\{\|x\|_p,\|y\|_p\})^2\qquad
\mbox{for all $x,y\in B_1^p(0)$}
\end{equation}
(cf.\ Lemma~\ref{secorder});
we may assume that $p\geq s$.
Then
$\|g(j(x,y))\|_q\leq
\|g(j(x,y))\|_r\leq \half (\|j(x,y)\|_s)^\alpha
\leq \half (\max\{\|x\|_p,\|y\|_p\})^{2\alpha}$
for all $x,y\in B_1^p(0)$,
by (\ref{usehld}) and (\ref{applty}).
Furthermore,
$\|g(x)\|_r\leq \half (\|x\|_s)^\alpha
\leq (\|x\|_s)^\alpha \leq (\|x\|_p)^\alpha$,
likewise
$\|g(y)\|_r\leq (\|y\|_p)^\alpha$,
and $\|g(j(x,y))\|_r
\leq \half (\max\{\|x\|_p,
\|y\|_p\})^{2\alpha}\leq (\max\{\|x\|_p,\|y\|_p\})^\alpha$,
entailing that
$\|D(g(x),g(y),g(j(x,y)))\|_q\leq \half (\max\{\|x\|_p,\|y\|_p\})^{2\alpha}$,
by (\ref{usetaylr}).
We now obtain
$\|g(j(x,y))+ D(g(x),g(y),g(j(x,y)))\|_q\leq
(\max\{\|x\|_p,\|y\|_p\})^{2\alpha}$ for
all $x,y\in B_1^p(0)$, using the triangle inequality.
\end{proof}
\begin{la}\label{La3}
There exists a continuous linear map
$\Lambda\!: L(G)\to L(H)$ such that
$\lambda(x)=\Lambda(x)$ for all
$x\in A$.
\end{la}
\begin{proof}
If we can show that
\begin{equation}\label{summ}
\lambda(x+y)=\lambda(x)+\lambda(y)
\qquad\mbox{for all $x,y\in A$ such that $x+y\in
A$,}
\end{equation}
then,
by \cite[Cor.\,A.2.27]{HaM},
the continuous map
$\lambda$
extends to
a continuous
homomorphism of groups $\Lambda\!: L(G)\to L(H)$.
Being a continuous homomorphism
between real topological vector spaces,
$\Lambda$
will be continuous linear.

To prove (\ref{summ}),
fix $x,y\in A$ such that $x+y\in A$.
There is $n_0\in \N$ such that $2^{-n}x\in Z$
and $2^{-n}y\in Z$ for all $n\geq n_0$.
For any such $n$, (\ref{taytau})
shows that
\begin{eqnarray*}
g(2^{-n}(x+y))&=&
g(2^{-n}x+2^{-n}y)\\
&=&g(2^{-n}x)\,*\, g(2^{-n}y)\, *\,
g((2^{-n}y)^{-1}\!*\!(2^{-n}x)^{-1}\!*\!(2^{-n}x+2^{-n}y))\\
&=&
g(2^{-n}x)\,*\, g(2^{-n}y)\, *\,
g(j(2^{-n}x,2^{-n}y))\\
&=&
g(2^{-n}x)\,+\,g(2^{-n}y)\,+\, r_n\,,
\end{eqnarray*}
where
$r_n:=
g(j(2^{-n}x,2^{-n}y))+D(g(2^{-n}x),g(2^{-n}y),
g(j(2^{-n}x,2^{-n}y)))$. Thus
\begin{equation}\label{letpass}
2^ng(2^{-n}(x+y))-
2^ng(2^{-n}x)-2^ng(2^{-n}y)= 2^n r_n\qquad \mbox{for all $n\geq n_0$.}
\end{equation}
Note that the left hand side of (\ref{letpass})
converges to $\lambda(x+y)-\lambda(x)-\lambda(y)$ as $n\to \infty$.
Hence
$\lambda(x+y)=\lambda(x)+\lambda(y)$ will hold if we can show
that $2^nr_n\to 0$ in $L(H)$ as $n\to \infty$.
To this end, given a continuous seminorm~$q$
on $L(H)$, let $p$
be as in Lemma~\ref{Latech}.
There is $n_1\geq n_0$ such that
$2^{-n}x,2^{-n}y\in B_1^p(0)$ for all $n\geq n_1$.
For any such~$n$, the cited lemma yields
$\|2^nr_n\|_q=2^n\|r_n\|_q\leq 2^n
(\max\{\|2^{-n}x\|_p,\|2^{-n}y\|_p\})^{2\alpha}
\leq \big(2^{-(2\alpha-1)}\big)^n\cdot (\max\{\|x\|_p,\|y\|_p\})^{2\alpha}$,
which tends to $0$ as $n\to \infty$.
Thus $2^nr_n\to 0$.
\end{proof}
\begin{la}\label{La4}
$g$ is totally differentiable at~$0$, with $g'(0)=\Lambda$.
\end{la}
\begin{proof}
Given a continuous seminorm $q$ on $L(H)$,
Lemma~\ref{LaNull}
provides a continuous seminorm $p$ on $L(G)$
such that $B_1^p(0)\sub A$
and (\ref{desired}) holds.
Choosing $n:=0$ and letting $m\to \infty$
in the first half of (\ref{auchbrauch}), we find that
\[
\left\|
\Lambda(x)-g(x)\right\|_q\leq
c\|x\|_p^{2\alpha}\quad\mbox{for all $x\in B_1^p(0)$,}
\]
where $c:=2^{2\alpha-1}\sum_{k=1}^\infty 2^{-(2\alpha-1)k}<\infty$.
Since $2\alpha-1>0$,
given $\ve>0$, there exists $\rho\in \;]0,1]$
such that $c \rho^{2\alpha-1}\leq \ve$.
Then $B_\rho^p(0)\sub A$, and for each $x\in B_\rho^p(0)$
we have
\[
\|g(x)-g(0)-\Lambda(x)\|_q
=\|g(x)-\Lambda(x)\|_q\leq c\|x\|^{2\alpha-1}_p\|x\|_p
\leq c\rho^{2\alpha-1}\|x\|_p\leq \ve \|x\|_p\,.
\]
Hence $g$ is totally differentiable at $0$, with
$g'(0)=\Lambda$. This completes the
proof of Lemma~\ref{La4}.
\end{proof}
Having proved Lemma~\ref{La4},
also Theorem~\ref{main}
is now fully established.
\end{proof}
Note that Lemma~\ref{impalph}
does not make use of the Mackey completeness of $L(H)$.
Beyond the real case (and independent
of Mackey completeness of $L(H)$), we still have:
\begin{prop}\label{globalH}
Let $\K$ be a valued field, $\alpha\in \,]0,1]$,
and $f\!: G\to H$ be a homomorphism
between Lie groups modelled on polynormed
$\K$-vector spaces.
Then $f$ is H\"{o}lder continuous
of degree $\alpha$ at~$1$ if and only if
$f$ is H\"{o}lder continuous of degree~$\alpha$.
\end{prop}
See Appendix~\ref{appb}
for the precise definitions and the proof.
\section{Homomorphisms between {\boldmath $p$\/}-adic
Lie groups}\label{secpad}
We now formulate a (slightly weaker)
analogue of Theorem~\ref{main} for $p$-adic
Lie groups. The proof carries over rather
directly,
whence we only indicate the most important changes.
\begin{thm}\label{varmain}
Let $f\!: G\to H$ be a homomorphism between
smooth Lie groups modelled
on polynormed $\Q_p$-vector spaces.
If $f$ is H\"{o}lder continuous at~$1$ and the modelling space
of~$H$ is Mackey complete,
then $f$ is $C^1_{\Q_p}$.
\end{thm}
{\bf Proof.} By hypothesis, $f$ is $H_\alpha$
at~$1$ for some $\alpha\in \;]0,1]$.
We let $\phi\!: U_1\to U\sub L(H)$
be a chart of~$H$ around~$1$, such that $\phi(1)=0$.
There exist open, symmetric
identity neighbourhoods $V_1\sub U_1$ and $W_1\sub V_1$
such that $V_1V_1\sub U_1$ and $(W_1)^{2p+1}:=
\underbrace{W_1W_1\cdots W_1}_{2p+1}\sub V_1$;\vspace{-5mm}
let $V:=\phi(V_1)$ and $W:=\phi(W_1)$. Define
\[
\mu\!: V\times V\to U,\qquad \mu(x,y):=x*y:=\phi(\phi^{-1}(x)\phi^{-1}(y))\,.
\]
Then the
$k$-fold products $x_1*x_2*\cdots * x_k$
are defined (and contained in~$V$),
for all $k\leq 2p+1$, $x_1,\ldots,x_k\in W$,
and every choice of brackets in this product.
The map
$\sigma\!: W\times W \to U$, $\sigma(x,y):=x^p*y$
satisfies $\sigma(0,0)=0$ and $\sigma'(0,0)(u,v)=pu+v$
for $u,v\in L(H)$.
The Taylor expansion around
$(0,0)$ yields
\[
\sigma(x,y)=px+y+R(x,y)\qquad\mbox{for all $x,y\in W$,}
\]
where $R(x,y):=R_1((0,0),(x,y),1)$.
Let $\psi\!: P_1\to P\sub L(G)$ be a chart of $G$
around~$1$, such that $f(P_1)\sub U_1$ and $\psi(1)=0$;
let $Q_1\sub P_1$ and $B_1\sub Q_1$ be symmetric
identity neighbourhoods such that $Q_1Q_1\sub P_1$,
$f(Q_1)\sub V_1$, $(B_1)^{2p+1}\sub Q_1$, and $f(B_1)\sub W_1$.
Set $Q:=\psi(Q_1)$ and $B:=\psi(B_1)$.
Define
$\nu\!: Q\times Q\to P$, $\nu(x,y):=x*y:=\psi(\psi^{-1}(x)\psi^{-1}(y))$.
Then
$g:=\phi\circ f|_{P_1}^{U_1}\circ \psi^{-1}\!: P \to U$
is continuous, maps $0$ to $0$,
and is $H_\alpha$ at~$0$.
Let
$A\sub B$ be a balanced, open
$0$-neighbourhood such that $g(x^{-p}*px)\in W$
for all $x\in A$. We define
\begin{equation}\label{defnh2}
h\!: A\to W,\quad h(x):=
g(x^{-p}* px)\, .
\end{equation}
For $x\in A$, we have
$g(px)=g(x^p*(x^{-p}*px))=g(x)^p*g(x^{-p}*px)
=\sigma(g(x),h(x))
=pg(x)+h(x)+R(g(x),h(x))$.
Likewise,
$g(p^2x)=pg(px)+h(px)+R(g(px),h(px))
=p^2g(x)+ph(x)+pR(g(x),h(x))+h(px)+R(g(px),h(px))$
and similarly
\begin{equation}\label{trick3}
g(p^n x)=
p^ng(x)+\sum_{k=1}^n
p^{n-k}
\left(
h(p^{k-1}x)+
R(g(p^{k-1}x),
h(p^{k-1}x))\right)
\end{equation}
for all $x\in A$ and $n\in \N_0$, by induction.
Hence
\begin{equation}\label{trick4}
\frac{g(p^n x)}{p^n}=
g(x)+\sum_{k=1}^n
p^{-k}
\left(
h(p^{k-1}x)+
R(g(p^{k-1}x),
h(p^{k-1}x))\right)
\end{equation}
for all $x\in A$ and
$n\in \N_0$. As in the proof of Lemma~\ref{LaNull},
we see:
\begin{la}\label{LaNull2}
Let $q$ be a continuous seminorm on $L(H)$.
Then there exists a continuous seminorm
$b$ on $L(G)$
such that $B_1^b(0) \sub A$,
\[
\|h(x)+R(g(x), h(x))\|_q\leq \|x\|_b^{2\alpha}\qquad
\mbox{for all $\, x\in B_1^b(0)$,}
\]
and $\|g(x)\|_q\leq (\|x\|_b)^\alpha$ for all $x\in B_1^b(0)$.\Punkt
\end{la}
Using Lemma~\ref{LaNull2},
we obtain by
a simple adaptation of the proof
of Lemma~\ref{impalph}
(where now $p\in \Q_p$
with $|p|=p^{-1}$ plays the role of $\half \in \R$):
\begin{la}\label{impalph2}
If $f$ is $H_\alpha$
at $1$ with
$\alpha\in \;]0,\frac{1}{2}]$,
then
$f$ also is
$H_{\frac{3}{2}\alpha}$ at~$1$.\Punkt
\end{la}
By the preceding,
we may assume now that
$\alpha\in \;]\half,1]$.
\begin{la}\label{La21}
The limit
$\lambda(x):=\lim_{n\to\infty} \frac{g(p^nx)}{p^n}$
exists in $L(H)$, for each $x\in A$.
For each continuous seminorm
$q$ on $L(H)$,
the convergence of
$\frac{g(p^nx)}{p^n}$
in $(L(H),\|.\|_q)$ is locally uniform
in~$x$. The map $\lambda\!: A\to L(H)$ is continuous.
\end{la}
\begin{proof}
The arguments from the real case are easily adapted.
To prove that $v_n:=p^{-n}g(p^nx_0)$
is a Mackey-Cauchy sequence for $x_0\in A$,
pick $0<\theta\in \Q$ such that $p^{-\theta} \in \, ]p^{-(2\alpha-1)},1[$;
set $r_{n,m}:=p^{[\theta(\smin\{n,m\}+1)]}\in \Q_p$,
where $[.]$ is the Gauss bracket (integer part).
Thus $|r_{n,m}|=p^{-[\theta(\smin\{n,m\}+1)]}\to 0$
as $n,m\to\infty$. Now complete the proof
as above.
\end{proof}
Let $Z\sub A$ be an open $0$-neighbourhood
such that $Z+Z\sub A$ and $Z^{-1}*Z^{-1}*(Z+Z)\sub A$.
We define
$j\!: Z\times Z\to A$,
$j(x,y):=y^{-1}*x^{-1}*(x+y)$.
Then $j(Z\times Z)\sub g^{-1}(W)$.
The map $\tau\!:
W\times W\times W\to U$, $\tau(x,y,z):=x*y*z$
is smooth, with $\tau(0,0,0)=0$ and
$\tau'(0,0,0)(u,v,w)=u+v+w$
for all $u,v,w\in L(H)$.
Let
$\wt{R}_1\!: (W\times W\times W)^{[1]}\to L(H)$
be the
first order Taylor remainder
of $\tau$.
Then
$\tau(x,y,z)=x+y+z + D(x,y,z)$
for all $x,y,z\in W$,
with $D(x,y,z):=\wt{R}_1((0,0,0),\, (x,y,z),\, 1)$.
Lemma~\ref{Latech} carries over:
\begin{la}\label{Latech22}
For every continuous seminorm $q$
on $L(H)$,
there is a continuous seminorm~$b$
on $L(G)$
such that $B_1^b(0) \sub Z$ and
\[
\left\|\,  g(j(x,y))\, + \, D \bigl( g(x),\,
g(y),\, g(j(x,y)) \bigr) \, \right\|_q \leq (\max\{\|x\|_b,
\|y\|_b\})^{2\alpha}\quad
\mbox{{\em for all $x,y\in B_1^b(0)$.\Punkt}}
\]
\end{la}
\begin{la}\label{La23}
$\lambda$ extends to a continuous $\Q_p$-linear map
$\Lambda\!: L(G)\to L(H)$.
\end{la}
\begin{proof}
The proof of Lemma~\ref{La3}
is easily adapted.
\end{proof}
In view of Lemma~\ref{totsmooth},
Theorem~\ref{varmain} follows
from the next lemma, whose proof directly parallels
that of Lemma~\ref{La4}:
\begin{la}\label{La24}
$g$ is totally differentiable at~$0$, with $g'(0)=\Lambda$.\Punkt
\end{la}
\appendix
\section{Proofs for the auxiliary results from
Section~\ref{secbase}}\label{secapp}
In this appendix, we prove
the results stated without proof in
Section~\ref{secbase}.
Not all techniques
from the real case carry over to general
valued fields~$\K$,
whence some of the proofs may look slightly unfamiliar.
In particular,
given an element $x$
of a polynormed $\K$-vector space~$E$
and a continuous seminorm $q$ on~$E$ such that
$\|x\|_q>0$, there need not be an element
$r\in \K$ such that $\|rx\|_q=1$.
As a substitute for normalization,
we shall frequently fix an element $a\in \K^\times$
such that $|a|<1$, and consider $a^{-k}x$
where $k\in \Z$ is chosen such that $|a|^{k+1}\leq \|x\|_q<|a|^k$.\\[3mm]
{\bf Proof of Lemma~\ref{secorder}.}
We use the second order Taylor expansion of~$f$,
\[
f(x+ty)-f(x)-tdf(x,y)=t^2 a_2(x,y)+t^2\,R_2(x,y,t)\quad
\mbox{for $(x,y,t)\in U^{[1]}$.}
\]
Since $R_2(x_0,0,0)=0$ and $a_2(x_0,0)=0$,
there exists $\rho\in \;]0,1]$
and a continuous seminorm $s$ on $E$
such that $B_{2\rho}^s(x_0)\sub U$,
\[
\|R_2(x,y,t)\|_q \leq 1\qquad
\mbox{for all $x\in B_\rho^s(x_0)$, $y\in B_\rho^s(0)$
and $|t|<\rho$,}
\]
and
$\|a_2(x,y)\|_q \leq 1$ for all $x\in B_\rho^s(x_0)$ and
$y\in B_\rho^s(0)$.
Pick $a\in \K^\times$ such that $|a|<1$;
define $\delta:=\rho^2|a|<\rho$,
$c:=2/(\rho|a|)^2$, and $p:= \max\{\frac{1}{\rho},\sqrt{\frac{c}{C}}\} s$.
Let $x\in B_1^p(x_0)$ and $y\in B_1^p(0)$;
then $x\in B_\delta^s(x_0)$ and $y\in B_\delta^s(0)$.
If $\|y\|_s >0$, there exists
$k\in \Z$ such that
$|a|^{k+1}\leq \rho^{-1}\|y\|_s  <|a|^k$.
Then $\|a^{-k}y\|_s<\rho$ and $|a^k|\leq |a|^{-1}\rho^{-1}\|y\|_s<\rho$.
If $\|y\|_s=0$, let $\ve\in \,]0,\rho[$
and choose $k\in \N$ so large that $|a|^k<\rho$ and
$2|a|^{2k}<\ve$.
Then, in either case,
\[
f(x+y)-f(x)=f(x+a^ka^{-k}y)-f(x)=df(x,y)
+a^{2k}a_2(x,a^{-k}y)+a^{2k}R_2(x,a^{-k}y,a^k)
\]
with
$r\!:=\!\|a^{2k}a_2(x,a^{-k}y)+a^{2k}R_2(x,a^{-k}y,a^k)\|_q \leq |a|^{2k}
(\|a_2(x,a^{-k}y)\|_q
+\|R_2(x,a^{-k}y,a^k)\|_q )$
$\leq 2|a|^{2k}$.
If $\|y\|_s>0$, the preceding formula shows that
$r\leq 2|a|^{-2}\rho^{-2}\|y\|^2_s
=c\|y\|^2_s\leq C\|y\|_p^2$.
If $\|y\|_s=0$, we have $r<\ve$ and thus $r=0\leq C\|y\|_p^2$,
as $\ve$ was arbitrary.
Hence $\|f(x+y)-f(x)-df(x,y)\|_q\leq C\|y\|_p^2$
for all $x\in B^p_1(x_0)$ and $y\in B_1^p(0)$.\vspace{3mm}\Punkt

\noindent
{\bf Proof of Lemma~\ref{basehold}.}
(a) and (b) are trivial;
(c) follows from Lemma~\ref{baseHglob}\,(a) and (e).\vspace{1mm}

(d) Let $E$, $F$ and $H$ be polynormed $\K$-vector spaces,
$U\sub E$ and $V\sub F$ be open, $x\in U$
and $f\!: U\to F$, $g\!: V\to H$ be maps such that
$f(U)\sub V$, $f$ is $H_\alpha$ at~$x$, and $g$ is $H_\beta$
at $f(x)$.
Given a continuous seminorm $q$ on~$H$, there exists
a continuous seminorm $p$ on $F$ such that
$B^p_1(f(x))\sub V$ and $\|g(z)-g(f(x))\|_q\leq \|z-f(x)\|_p^\beta$
for all $z\in B^p_1(f(x))$.
There is a continuous
seminorm $r$ on $E$ such that $B^r_1(0)\sub U$
and $\|f(y)-f(x)\|_p\leq \|y-x\|_r^\alpha\leq 1$
for all $y\in B^r_1(x)$.
Then $\|g(f(y))-g(f(x))\|_q\leq \|f(y)-f(x)\|_p^\beta
\leq \|y-x\|_r^{\alpha \beta}$
for all $y\in B^r_1(0)$.\vspace{3mm}\Punkt

\noindent
{\bf Proof of \ref{claim1tot}.}
{\em If $h$ is tangent to $0$,}
let $q$ be a continuous seminorm on~$F$.
For $W:=B_1^q(0)$ we then find $V$ and $\theta\!: I\to \K$
as in {\bf \ref{dftot}}. We may assume
that $V$ is balanced and $I=B_r(0)\sub \K$
for some $r>0$.
There exists a continuous seminorm $p$ on~$E$
such that $B_1^p(0)\sub V$
and $B_1^p(0)\sub U-x$.
Replacing $V$ with $B_1^p(0)$,
we may assume that $V=B^p_1(0)$.
Fix $a\in \K^\times$ such that $|a|<1$.
Given $\ve>0$, there exists $\delta\in \,]0,1]$
such that $\frac{|\theta(t)|}{|t|}<\ve|a|$ if $|t|<\delta$.
Then $B^p_\delta(0)\sub B^p_1(0)\sub U-x$.
Let $y\in B^p_\delta(0)$;
we claim that $\|h(y)\|_q\leq \ve\|y\|_p$.
If $\|y\|_p=0$, then $t^{-1}y\in V$ for each
$0\not=t\in I$, whence $h(y)=h(t(t^{-1}y))\in \theta(t)W$
and thus $\|h(y)\|_q\leq |\theta(t)|$. Hence
$\|h(y)\|_q=0\leq \ve\|y\|_p$.
If $\|y\|_p>0$,
then there is $k\in \N_0$
such that $|a|^{k+1}\leq \|y\|_p <|a|^k$.
Set $t:=a^k$. Then $t^{-1} y\in V$ and
thus $h(y)=h(t(t^{-1}y))\in \theta(t)W$,
whence $h(y)=\theta(t)w$ with $w\in W$.
Hence $\|h(y)\|_q=|\theta(t)|\cdot
\|w\|_q\leq |\theta(t)|\leq \ve \,|a| \,|t| \leq \ve \|y\|_p$.\vspace{1mm}

{\em Conversely},
assume that the condition from {\bf \ref{claim1tot}}
is satisfied. Given a $0$-neighbourhood
$W\sub F$,
there exists a continuous seminorm $q$ on $F$ such that
$B_1^q(0)\sub W$.
We choose a continuous seminorm $p$ on~$E$ as
described in {\bf \ref{claim1tot}}.
Let
$(a_n)_{n\in \N}$ be a sequence in $\K^\times$ such that
$\lim_{n\to\infty} a_n=0$.
For each~$n$, there exists
$\delta_n>0$ such that $B_{\delta_n}^p(0)\sub U-x$
and $\|h(y)\|_q\leq |a_n|\cdot \|y\|_p$
for all $y\in B_{\delta_n}^p(0)$.
We may assume that $\delta_1>\delta_2>\cdots$ and $\lim_{n\to\infty}
\delta_n=0$. Now set $V:=B_1^p(0)$ and
define $\theta\!: I\to\K$
on $I:=B_{\delta_1}(0)\sub \K$ via $\theta(0)=0$,
$\theta(t):=a_nt$ if $|t|\in [\delta_{n+1},\delta_n[$.
Then $\theta(t)=o(t)$ and $IV\sub B_{\delta_1}^p(0)\sub U-x$.
Furthermore, $h(tv)\sub \theta(t)W$ for
$t\in I$ and $v\in V$:
This is trivial if $t=0$, and also
if $t\not=0$ and $\|v\|_p=0$,
because then $\|h(tv)\|_q\leq |a_n|\cdot \|tv\|_p=0$
(for any~$n$)
and hence $h(tv)\in q^{-1}(\{0\})\sub \theta(t)W$.
Otherwise, $0<\|tv\|_p\in [\delta_{n+1},\delta_n[$
for some $n$ and thus $\|h(tv)\|_q\leq |a_n|\cdot \|tv\|_q
= |t a_n|\cdot \|v\|_q<|ta_n|$,
whence $h(tv)\in B_{|t a_n|}^q(0)=ta_n B_1^q(0)=\theta(t)B_1^q(0)\sub
\theta(t)W$.
Hence $h$ is tangent to~$0$.\vspace{3mm}\Punkt

\noindent
{\bf Proof of \ref{claim2tot}.}
Let $f$, $g$ be as in {\bf \ref{claim2tot}},
and $W_1\sub H$ be a $0$-neighbourhood.
There is a balanced $0$-neighbourhood $W\sub H$
such that $W+W\sub W_1$.
As $g$ is totally differentiable at~$f(x)$
and $g'(f(x))$ continuous linear,
we find balanced $0$-neighbourhoods $P_1 \sub F$,
$I\sub \K$ and a map $\theta\!: I\to\K$
which is $o(t)$, such that
$g'(f(x)).P_1\sub W$,
$I P_1\sub V-f(x)$,
and $h_2(t P_1)\sub \theta(t)W$ for $t\!\in\! I$,
where $h_2\!: \!V\!\!-\!f(x)\to H$,
$h_2(z)\!=g(f(x)+z)-g(f(x))-g'(f(x)).z$.
There is a balanced $0$-neighbourhood $P\sub F$
such that $P+P\sub P_1$.

Define $h_1\!: U-x\to F$, $h_1(y):=f(x+y)-f(x)-f'(x).y$.
There are $0$-neighbourhoods $Q\sub E$,
$J\sub \K$ and a map $\xi\!: J\to \K$
which is $o(t)$,
such that $JQ\sub U-x$, $f'(x).Q\sub P$,
and $h_1(tQ)\sub \xi(t)P$.
After shrinking $I$ and $J$, we may assume
that $I=J$ and
$\big|\frac{\xi(t)}{t}\big|\leq 1$
for all $0\not=t\in I$.
Define $\eta\!: I\to \K$ via
$\eta(t):=\theta(t)$ if $|\theta(t)|\geq |\xi(t)|$,
$\eta(t):=\xi(t)$ if $|\theta(t)|<|\xi(t)|$.
Define
$A:=g'(f(x))\circ f'(x)$ and $h\!: U-x\to H$,
$h(y):=g(f(x+y))-g(f(x))-A.y$. Then
\begin{eqnarray*}
h(y) &=& g(f(x)+ f'(x).y+h_1(y))\, - g(f(x))-A.y\\
&=& g(f(x))+g'(f(x)).z +h_2(z)-g(f(x))-A.y\\
&=& g'(f(x)).h_1(y)+ h_2(f'(x).y+h_1(y))\,,
\end{eqnarray*}
where $z:=f'(x).y+h_1(y)$.
Let $t\in I$ and $y\in Q$.
Then $h_1(ty)\in \xi(t)P\sub \xi(t)P_1\sub \eta(t)P_1$
as $P_1$ is balanced,
and thus $g'(f(x)).h_1(ty)\in \eta(t)W$.
Furthermore,
$f'(x).ty\in t P$ and $h_1(ty)\in \xi(t) P\sub
t P$ (as $|\xi(t)|\leq |t|$),
whence $f'(x).ty+h_1(ty)\in t(P+ P)\sub tP_1$
and thus $h_2(f'(x).ty+h_1(ty))\in \theta(t)W\sub \eta(t)W$,
using that $W$ is balanced.
Hence $h(ty)
= g'(f(x)).h_1(ty)+ h_2(f'(x).ty+h_1(ty))\in \eta(t)(W+W)\sub \eta(t)W_1$,
and thus $h(tQ)\sub \eta(t)W_1$.
We have shown that $h$ is tangent to~$0$;
the assertions follow.\vspace{3mm}\Punkt

\noindent
{\bf Proof of Lemma~\ref{C2total}.}
We consider the second order Taylor expansion of $f\,$:
\begin{equation}\label{use2tay}
f(x+tv)=f(x)+t\, df(x,v) +t^2 a_2(x,v)+t^2 R_2(x,v,t)\qquad \mbox{for
all $(x,v,t)\in U^{[1]}$}
\end{equation}
(see Proposition~\ref{basictaylor}).
Fix $x\in U$.
The map $f'(x):=df(x,\sbull)\!: E\to F$
being continuous linear,
to establish total differentiability of~$f$ at $x$ we
only need to show that
\[
h\!: U-x\to F, \quad h(y):=f(x+y)-f(x)-f'(x).y
\]
is tangent to~$0$.
To this end, let $W$ be a $0$-neighbourhood in~$F$.
There exists a $0$-neighbourhood $W_1\sub F$
such that $W_1+W_1\sub W$.
As $R_2(x,0,0)=0$
and $R_2$ is continuous,
there is a $0$-neighbourhood $V\sub E$
and a $0$-neighbourhood $I\sub \K$
such that $(x,v,t)\in U^{[1]}$
and $R_2(x,v,t)\in W_1$
for all $v\in V$ and $t\in I$.
Since $a_2$ is continuous
and $a_2(x,0)=0$, after shrinking $V$ we may assume
that furthermore $a_2(x,v)\in W_1$ for all $v\in V$.
Define
\[
\theta\!: I\to \K,\quad \theta(t):=t^2\,.
\]
Then $\theta(t)=o(t)$.
For each $t\in I$
and $y\in tV$, say $y=tv$ with $v\in V$,
we have
\[
h(y)=t^2(a_2(x,v)+R_2(x,v,t))
\in t^2(W_1+W_1)\sub t^2W=\theta(t)W\,,
\]
using (\ref{use2tay}).
Hence $h$ is indeed tangent to~$0$.\vspace{3mm}\Punkt

\noindent
{\bf Proof of Lemma~\ref{2variable}.}
If $f$ is $C^1$, set $f'(x):=df(x,\sbull)$. Then
$\wt{f}(y,t)=f^{[1]}(x,y,t)$
is continuous.\vspace{1mm}

Now assume that $\K$ is a valued field,
$f$ is continuous on~$U$
and totally differentiable at~$x$.
Define $\wt{f}_x\!: \wt{U}_x\to E$
as in {\bf \ref{defnfeeb}},
using the total differential $f'(x)$.
Since $f$ is continuous,
so is $\wt{f}_x|_A$.
By a theorem
of Bourbaki and Dieudonn\'{e}
\cite[Exerc.\,3.2\,A\,(b)]{Eng}, the map $\wt{f}_x$ is continuous
if its restriction
$\wt{f}_x|_{A\cup \{(y,0)\}}$
is continuous for each $y\in E$.
This will hold if we can show that $\wt{f}_x(y_\alpha,t_\alpha)
\to \wt{f}_x(y,0)$ for each net $(y_\alpha,t_\alpha)$
in $A$ converging to $(y,0)$ for some $y\in E$.
To see that this condition is satisfied,
let $W_1\sub F$ be a $0$-neighbourhood.
There is a balanced $0$-neighbourhood $W\sub F$
such that $W+W\sub W_1$.
Since $f$ is totally differentiable at~$x$,
there exists an open $0$-neighbourhood $V\sub E$
and a function $\theta\!: I\to \K$
on some $0$-neighbourhood in~$\K$ such that $I\cdot V\sub U-x$ holds,
$\theta(t)=o(t)$, and
\begin{equation}\label{usetotal}
f(x+sv)\in f(x)+s f'(x).v + \theta(s)W\qquad
\mbox{for all $v\in V$ and $s\in I$.}
\end{equation}
Pick $r\in \K^\times$ such that $r y\in V$.
As $(y_\alpha,t_\alpha)\to (y,0)$,
there exists $\beta$ such that $f'(x).(y_\alpha-y)\in W$,
$v_\alpha:=r y_\alpha\in V$,
$s_\alpha:=r^{-1}t_\alpha\in I$,
and $|\theta(s_\alpha)|/|s_\alpha|\leq |r|$
for all $\alpha\geq \beta$.
For any such $\alpha$,
(\ref{usetotal})
applied to $x+t_\alpha y_\alpha=x+s_\alpha v_\alpha$
shows that
\[
\wt{f}_x(y_\alpha,t_\alpha)-
\wt{f}_x(y,0)\in f'(x).y_\alpha-f'(x).y+\frac{\theta(s_\alpha)}{t_\alpha}W
\sub W+ \frac{\theta(s_\alpha)}{r s_\alpha}W
\sub W+W\sub W_1\,.
\]
Thus indeed $\wt{f}_x(y_\alpha,t_\alpha)\to
\wt{f}_x(y,0)$.\vspace{2mm}\Punkt

\noindent
{\bf Proof of \ref{chainfeeb}.}
We define $\wt{f}_x\!: \wt{U}_x\to F$ and $\wt{g}_{f(x)}\!: \wt{V}_{f(x)}\to H$
as in {\bf \ref{defnfeeb}} and
abbreviate $h:=g\circ f\!: U\to H$.
For any $y\in E$ and $t\in \K^\times$ such that
$x+ty\in U$, we calculate
\[
\frac{h(x+ty)-h(x)}{t} \,=\,
\frac{g\Big( f(x)+t\, \frac{f(x+ty)-f(x)}{t}\Big)-g(f(x))}{t}
\,=\, \wt{g}_{f(x)}\big(\wt{f}_x(y,t),t\big)\,=\, \wt{h}_x(y,t)
\]
where $\wt{h}_x\!: \wt{U}_x\to H$,
$\wt{h}_x(y,t):=\wt{g}_{f(x)}\big(\wt{f}_x(y,t),t\big)$
is continuous, and the map $\wt{h}_x(\sbull,0)=g'(f(x))\circ f'(x)$
is continuous linear. Thus $h=g\circ f$ is feebly
differentiable at~$x$.\Punkt
\section{H\"{o}lder continuity at 1 entails H\"{o}lder
continuity}\label{appb}
So far, we only considered H\"{o}lder continuity
{\em at a point.} We now discuss
mappings which are H\"{o}lder continuous
on all of their domain.
Basic facts are provided
and a proof for Proposition~\ref{globalH}
is given.
\begin{defn}
Let $E$ and $F$
be polynormed vector spaces over a valued field~$\K$,
and $U\sub E$ be open.
A map $f\!: U\to F$
is called {\em H\"{o}lder continuous of degree $\alpha$\/}
(or $H_\alpha$, for short)
if, for every $x_0\in U$ and continuous seminorm $q$
on $F$, there exists a continuous seminorm $p$
on $E$ and $\delta>0$ such that $B_\delta^p(x_0)\sub U$ and
$\|f(y)-f(x)\|_q\leq \|y-x\|_p^\alpha$,
for all $x,y\in B_\delta^p(x_0)$.
If $f$ is $H_1$, we also say that $f$ is {\em Lipschitz continuous\/}.
\end{defn}
\begin{la}\label{baseHglob}
For maps between open subsets of polynormed
$\K$-vector spaces, we have:
\begin{itemize}
\item[\rm (a)]
If $f\!: E\supseteq U\to F$ is $H_\alpha$,
then $f$ is $H_\alpha$ at each $x\in U$.
\item[\rm (b)]
If $f$ is $H_\alpha$ then $f$ is continuous.
\item[\rm (c)]
If $\alpha\geq \beta$ and $f$ is $H_\alpha$, then
$f$ is $H_\beta$.
\item[\rm (d)]
If $f$ and $g$ are composable maps
such that $f$ is $H_\alpha$
and $g$ is $H_\beta$,
then $g\circ f$ is $H_{\alpha\cdot \beta}$.
\item[\rm (e)]
Any $C^1$-map is Lipschitz continuous.
\end{itemize}
\end{la}
\begin{proof}
(a), (b) and (c) are obvious;
(d) can be proved as Lemma~\ref{basehold}\,(d).\vspace{1mm}

(e) We use the first order Taylor expansion
$f(x+ty)-f(x)=tdf(x,y)+t R_1(x,y,t)$ of the $C^1$-map
$f\!: E\supseteq U\to F$. Here
\[
R_1(x,y,1)=tR_1(x,t^{-1}y,t)\quad \mbox{for $t\in \K^\times$
and $(x,y)\in U\times E$
such that $x+y\in U$.}
\]
Fix $x_0\in U$. Let $q$ be a continuous seminorm on~$F$.
Pick $a\in \K^\times$ such that $|a|<1$.
Since $df(x_0,0)=0$, using the continuity of $df$
we find a continuous seminorm $r$ on~$E$
such that $B_1^r(x_0)\sub U$ and
$\|df(x,y)\|_q\leq |a|$
for all $x\in B^r_1(x_0)$ and $y\in B^r_1(0)$,
whence $\|df(x,y)\|_q\leq \|y\|_r$
for all $x\in B^r_1(x_0)$ and $y\in E$.
Since $R_1(x_0,0,0)=0$, we find a continuous
seminorm $p$ on~$E$ and $\rho\in \,]0,1]$ such that
$B^p_{2 \rho}(x_0)\sub U$ and
$\|R_1(x,y,t)\|_q\leq 1$
for all $x\in B^p_\rho(x_0)$,
$y\in B^p_\rho(0)$ and $t\in B_\rho(0)\sub \K$;
we may assume that $p\geq r$.
Define $\delta:=\half \rho^2|a|$.
Given $x,y\in B^p_\delta(x_0)$, set $z:=y-x$.
If $\|z\|_p>0$, there is $k\in \Z$ such that
$|a|^{k+1}\leq \rho^{-1}\|z\|_p<|a|^k$.
Then $\|a^{-k}z\|_p<\rho$ and $|a^k|\leq |a|^{-1}\rho^{-1}\|z\|_p<\rho$,
whence $\|R_1(x,z,1)\|_q=|a^k|\, \|R_1(x,a^{-k}z,a^k)\|_q
\leq |a^k|\leq |a|^{-1}\rho^{-1}\|z\|_p$ and thus
$\|f(x+z)-f(x)\|_q\leq \|df(x,z)\|_q+\|R_1(x,z,1)\|_q
\leq (1+|a|^{-1}\rho^{-1})\|z\|_p$. Hence
\begin{equation}\label{fn2}
\|f(y)-f(x)\|_q \leq (1+|a|^{-1}\rho^{-1})\|y-x\|_p\,.
\end{equation}
If $\|z\|_p=0$, given $\ve>0$ pick $t\in \K^\times$
such that $|t|<\min\{\rho,\ve\}$.
Then $\|df(x,z)\|_q=0$
and $\|R_1(x,z,1)\|_q=|t|\, \| R_1(x,t^{-1}z,t)\|_q\leq |t|\leq \ve$,
whence $\|R_1(x,z,1)\|_q=0$ (as $\ve$ was arbitrary).
Thus (\ref{fn2}) also holds if $\|z\|_p=0$.
\end{proof}
\begin{defn}
Let $f\!: M\to N$ be a map between $C^1_\K$-manifolds
modelled on polynormed $\K$-vector spaces,
and $\alpha\in \;]0,1]$.
We say that $f$ is {\em H\"{o}lder continuous of degree
$\alpha$\/}
(or briefly: $f$ is $H_\alpha)$,
if $f$ is continuous
and, for each $x_0\in M$, there exist
a chart
$\phi\!: U_1\to U$ of $M$ around~$x_0$
and a chart $\psi\!: V_1\to V$ of $N$ around
$f(x_0)$, such that
$\,\phi(f^{-1}(V_1)\cap U_1)\to
V$,
$y\mto \psi(f(\phi^{-1}(y)))\,$
is $H_\alpha$.
\,(This then holds for any choice of
$\phi$ and $\psi$, by Lemma~\ref{baseHglob}).
\end{defn}
{\bf Proof of Proposition~\ref{globalH}.}
Let $f$ be $H_\alpha$ at~$1$.
If we can show that $f|_U$ is $H_\alpha$
for an open identity neighbourhood $U\sub G$,
then $f|_{xU}=\lambda^H_{f(x)}\circ f|_U\circ \lambda^G_{x^{-1}}|_{xU}^U$
(with left translation maps as indicated)
will be $H_\alpha$ by Lemma~\ref{baseHglob}\,(d) and (e),
for each $x\!\in \! G$, whence $f$ will~be~$H_\alpha$.\vspace{1mm}

We choose charts $\phi\!: P_1\to P\sub L(G)$
and $\psi\!: Q_1\to Q\sub L(H)$ around~$1$ of $G$ and $H$, respectively,
such that $\phi(1)=0$, $\psi(1)=0$ and $f(P_1)\sub Q_1$.
There are symmetric identity neighbourhoods
$X_1,U_1\sub G$ and $Y_1,V_1\sub H$ such that $X_1X_1\sub P_1$,
$U_1U_1\sub X_1$,
$Y_1Y_1\sub Q_1$, $V_1V_1\sub Y_1$,
$f(X_1)\sub Y_1$, and $f(U_1)\sub V_1$;
set $X:=\phi(X_1)$, $U:=\phi(U_1)$, $Y:=\psi(Y_1)$,
and $V:=\psi(V_1)$.
We write $\mu\!: X\times X\to P$ and $\nu\!: Y\times Y\to Q$
(or ``$*$'')
for the local multiplications obtained from the
respective group multiplication,
and define $g:=\psi \circ f|_P \circ \phi^{-1}|_P\!: P\to Q$.
Let $q$ be a continuous seminorm on $L(H)$,
and $x_0\in U$.
As~$\nu$ is $H_1$,
there is a continuous seminorm $p$
on $L(H)$ such that
$B^p_1(g(x_0))\times B^p_1(0)\sub V\times V$~and
\[
\|\nu(u,v)-\nu(u',v')\|_q\leq \max\{\|u'-u\|_p,\|v'-v\|_p\}\quad
\mbox{for all $u,u'\in B^p_1(x_0)$, $v,v'\in B^p_1(0)$.}
\]
Now $g$ being $H_\alpha$ at~$0$,
there exists a continuous seminorm $r$ on~$L(G)$
such that $B^r_1(0)\sub X$ and $\|g(x)\|_p\leq \|x\|^\alpha_r$
for all $x\in B^r_1(0)$.
The map $h\!: U\times U\to X$, $h(u,v):=u^{-1}*v$
being Lipschitz continuous,
there is a continuous seminorm $s\geq r$ on $L(G)$
such that $B^s_1(x_0)\sub U$,
$h(B^s_1(x_0)\times B^s_1(x_0))\sub B^r_1(0)$,
and
\[
\|h(u,v)-h(u',v')\|_r\leq \max\{\|u-u'\|_s,\|v-v'\|_s\}\quad
\mbox{for all $u,u',v,v'\in B^s_1(x_0)$.}
\]
For any $x,y\in B^s_1(x_0)\sub U$, we obtain
\begin{eqnarray*}
\|g(y)-g(x)\|_q & = & \|g(x)*g(x^{-1}*y)-g(x)\|_q=
\|\nu(g(x),g(x^{-1}*y))-\nu(g(x),0)\|_q\\
&\leq & \max\{\|g(x)-g(x)\|_p,\|g(x^{-1}*y)\|_p\}
=\|g(x^{-1}*y)\|_p\\
&\leq & \|x^{-1}*y\|_r^\alpha =\|h(x,y)-h(x,x)\|_r^\alpha
\leq \|y-x\|_s ^\alpha\! .
\end{eqnarray*}
Hence $g|_U$ is $H_\alpha$ indeed.\Punkt
\noindent
{\footnotesize
{\bf Helge Gl\"{o}ckner}, TU~Darmstadt, FB~Mathematik~AG~5,
Schlossgartenstr.\,7, 64289 Darmstadt, Germany.\\
E-Mail: gloeckner@mathematik.tu-darmstadt.de}

\begin{thebibliography}{10}
%
%
\bibitem{AaS} Averbukh, V.\,I. and O.\,G. Smolyanov, {\em The various
definitions of the derivative in linear topological spaces}, Russ.\ Math.\
Surv.\ {\bf 23}\,(1968), 67--113.
%
%
\bibitem{BGN} Bertram, W., H. Gl\"{o}ckner and K.-H. Neeb,
{\em Differential calculus over general base fields and rings},
to appear in Expo.\ Math.; also arXiv:math.GM/0303300\,.
%
%
\bibitem{Eng} Engelking, R., ``General Topology,''
Heldermann Verlag, 1989.
%
%
\bibitem{RES} Gl\"{o}ckner, H.,
{\em Infinite-dimensional Lie groups without
completeness restrictions},
pp.\,43--59 in: Strasburger, A. et al.\ (Eds.),
Geometry and Analysis on Finite- and Infinite-Dimensional
Lie Groups, Banach Center Publications Vol.~{\bf 55},
Warsaw, 2002.
%
%
\bibitem{INF} -----, ``Infinite-Dimensional Analysis,''
\LaTeX-lecture notes for a course held at Darmstadt University of
Technology (Winter Semester 2002--2003).
%
%
\bibitem{ANA} -----, {\em Every smooth $p$-adic Lie group
admits a compatible analytic structure}, TU Darmstadt Preprint {\bf 2307},
December 2003; also arXiv:math.GR/0312113\,.
%
%
\bibitem{FUN} -----, {\em Fundamentals of direct limit Lie theory},
TU Darmstadt Preprint {\bf 2324}, March 2004; also arXiv:math.GR/0403093\,.
%
%
\bibitem{CON} -----,
{\em Conveniently H\"{o}lder homomorphisms are
smooth in the convenient sense},
TU Darmstadt Preprint {\bf 2330}, April 2004;
also arXiv:math.GR/0404344.
%
%
\bibitem{ZOO} -----, {\em Lie groups over non-discrete
topological fields}, in preparation.
%
%
\bibitem{GaN} Gl\"{o}ckner, H. and K.-H. Neeb, {\em Banach-Lie quotients,
enlargibility, and universal complexifications},
J. Reine Angew.\ Math.\ {\bf 560}\,(2003), 1--28.
%
%
\bibitem{Ham} Hamilton, R., {\em The inverse function theorem of Nash
and Moser}, Bull.\ Amer.\ Math.\ Soc.\ {\bf 7}\,(1982), 65--222.
%
%
\bibitem{HaM} Hofmann, K.\,H. and S.\,A. Morris,
``The Structure of Compact Groups,''
de Gruyter, 
1998.
%
%
\bibitem{Kel} Keller, H.\,H., ``Differential Calculus
in Locally Convex Spaces,'' Springer-Verlag, 1974.
%
%
\bibitem{KaM} Kriegl, A. and P.\,W. Michor,
``The Convenient Setting of Global Analysis,''
Math.\ Surveys and Monographs {\bf 53},
AMS, Providence, 1997.
%
%
\bibitem{Lan} Lang, S. ``Fundamentals of Differential Geometry,''
Springer-Verlag, 
1999.
%
%
\bibitem{Mic} Michor, P.\,W., ``Manifolds of Differentiable Mappings,''
Shiva, 1980.
%
%
\bibitem{MiP} Milnor, J., {\em On infinite dimensional Lie groups},
Preprint, Institute for Advanced Study, Princeton, 1982.
%
%
\bibitem{Mil} -----, {\em Remarks on infinite-dimensional
Lie groups}, pp.\,1008--1057 in:
B. DeWitt and R. Stora (Eds.),
``Relativity, Groups and Topology~II,'' North-Holland, 1983.
%
%
\bibitem{Omo} Omori, H., {\em Groups of diffeomorphisms and their
subgroups}, Trans.\ Amer.\ Math.\ Soc.\ {\bf 178}\,(1973),
85--122.
%
%
\bibitem{Rob} Robart, T., {\em Sur l'int\'{e}grabilit\'{e}
des sous-alg\`{e}bres de Lie en dimension infinie},
Canad.\ J. Math.\ {\bf 49}\,(1997), 820--839.
%
%
\bibitem{Wie} Wi\c{e}s\l{}aw, W.,
``Topological Fields,'' Marcel Dekker,
1988.
%
%
\end{thebibliography}
\end{document}